\documentclass[12pt,twoside,leqno]{amsart}

\usepackage{amssymb,amsbsy,amsmath,amsfonts,eucal,amsthm,
amssymb,graphicx,times,xypic,color,wasysym,mathrsfs}

\usepackage[T1]{fontenc} 
\definecolor{blue}{cmyk}{1.,1.,0.,0.53}
\definecolor{red}{cmyk}{0.,1.,1.,0.53}
\definecolor{green}{cmyk}{1.,0.,1.,0.53}

\let\mathcal\mathscr

\newcommand{\C}{\mathbb{C}}\newcommand{\K}{\mathbb{K}}
\newcommand{\N}{\mathbb{N}}
\newcommand{\R}{\mathbb{R}}

\begin{document}


\title[Nonrigid spherical real analytic hypersurfaces in $\C^2$]{
Nonrigid spherical 
\\
real analytic hypersurfaces in $\C^2$
}

\author{Jo\"el Merker}

\address{D\'epartement de math\'ematiques 
et applications, \'Ecole Normale Sup\'erieure, Paris}
\email{merker@dma.ens.fr} 

\date{\number\year-\number\month-\number\day.\ To appear
in Complex Variables and Elliptic Equations}

\begin{abstract}
A Levi nondegenerate real analytic hypersurface $M$ of $\C^2$
represented in local coordinates $(z, w) \in \C^2$ by a complex
defining equation of the form $w = \Theta ( z, \overline{ z},
\overline{ w})$ which satisfies an appropriate reality condition, is
spherical if and only if its complex graphing function $\Theta$
satisfies an explicitly written sixth-order polynomial complex partial
differential equation. In the rigid case (known before), this system
simplifies considerably, but in the general nonrigid case, its
combinatorial complexity shows well why the two fundamental curvature
tensors constructed by \'Elie Cartan in 1932 in his classification of
hypersurfaces have, since then, never been reached in parametric
representation.
\end{abstract}

\maketitle

\begin{center}
\begin{minipage}[t]{11cm}
\baselineskip =0.35cm
{\scriptsize

\centerline{\bf Table of contents}

\smallskip

{\bf 1.~Introduction\dotfill 
\pageref{Section-1}.}

{\bf 2.~Segre varieties and differential equations\dotfill 
\pageref{Section-2}.}

{\bf 3.~Geometry of associated submanifolds of solutions\dotfill 
\pageref{Section-3}.}

{\bf 4.~Effective differential characterization of sphericality
in $\C^2$\dotfill 
\pageref{Section-4}.}

{\bf 5.~Some complete expansions: examples of 
expression swellings\dotfill 
\pageref{Section-5}.}

}\end{minipage}
\end{center}

\section*{\S1.~Introduction}
\label{Section-1}

A real analytic hypersurface $M$ in $\C^2$ is called {\sl spherical}
at one of its points $p$ if there exists a nonempty open neighborhood
$U_p$ of $p$ in $\C^2$ such that $M \cap U_p$ is biholomorphic to a
piece of the unit sphere $S^3 = \big\{ (z, w) \colon \vert z \vert^2 +
\vert w \vert^2 = 1 \big\}$. When $M$ is connected, sphericality at
one point is known to propagate all over $M$, for it is equivalent to
the vanishing of two certain {\em real analytic} curvature tensors
that were constructed by \'Elie Cartan in~\cite{ ca1932}. However, the
intrinsic computational complexity, in the Cauchy-Riemann (CR for
short) context, of \'Elie Cartan's algorithm to derive an absolute
parallelism on some suitable eight-dimensional principal bundle
$\mathcal{ P} \to M$ prevents from controlling explicitly all the
appearing differential forms. As a matter of fact, the effective
computation, in terms of a defining equation for $M$, of the two
fundamental differential invariants the vanishing of which
characterizes sphericality, appears nowhere in the literature ({\em
see}~ {\em e.g.}~\cite{ ns2003, crsa2005, isa2009} and the references
therein as well), except notably when one makes the assumption that,
in some suitable local holomorphic coordinates $(z, w) = ( x + iy, \,
u + iv)$ vanishing at the point $p$, the defining equation is of the
so-called {\sl rigid} form $u = \varphi ( x, y)$ with the variable $v$
missing, or even of the so-called (simpler) {\sl tube} form $u =
\varphi ( x)$, with the two variables $y$ and $v$ missing, {\em
see}~\cite{ isa2009} which showed recently a renewed interest, in CR
geometry, for explicit characterizations of sphericality.  But in
general, a real analytic hypersurface $M \subset \C^2$ is represented
at $p$ by a {\em real} equation $u = \varphi (x, y, v)$ whose graphing
function $\varphi$ depends entirely arbitrarily upon $v$ {\em also},
and then apparently, the characterization of sphericality is still
unknown.

On the other hand, in the studies~\cite{ me2001a, me2001b, me2002,
me2005a, me2005b} devoted to the CR reflection principle, it was
emphasized that all the adequate invariants of CR mappings between CR
manifolds: {\sl Pair of Segre foliations}, {\sl Segre chains}, {\sl
Complexified CR orbits}, {\sl Jets of complexified Segre varietes},
{\sl Rigidity of formal CR mappings}, {\sl Nondegeneracy conditions},
{\sl CR-reflection function}\footnote{\,
For a presentation of these
concepts, the reader is referred to the extensive introductions
of~\cite{ me2002, me2005b} and also to~\cite{ me2009} for more about
why dealing only with complex defining equations is natural and
unavoidable when one wants to insert CR geometry in the wider
universe of completely integrable systems of real or complex analytic
partial differential equations. },
can be viewed correctly only when $M$ is represented by a so-called
{\sl complex defining equation} of the form:
\[
w
=
\Theta\big(z,\overline{z},\overline{w}\big),
\]
where the function $\Theta \in \C \big\{ z, \overline{ z}, \overline{
w} \big\}$, vanishing at the origin, is the unique function obtained by
solving with respect to $w$ the equation: $\frac{ w + \overline{
w}}{2} = \varphi \big( \frac{z + \overline{ z} }{ 2},\, \frac{z -
\overline{ z} }{ 2i},\, \frac{w - \overline{ w} }{ 2i} \big)$; then
the fact that $\varphi$ was {\em real} is reflected, in terms of this
new function $\Theta ( z, \overline{ z}, \overline{ w})$, by the
constraint that, together with its complex conjugate $\overline{
\Theta} \big( z, \overline{ z}, \overline{ w} \big)$, it satisfies the
functional equation\footnote{\,
More will be said shortly in Section~2 below. }:
\[
w
\equiv
\Theta\big(z,\,\overline{z},\,
\overline{\Theta}(\overline{z},z,w)\big).
\]
Accordingly, the author suspected since a
few years\,\,---\,\,{\em cf.}~the
Open Question~2.35 in~\cite{ me2009}\,\,---\,\,that sphericality of
$M$ at $p$ should and could be expressed adequately in terms of
$\Theta$. The classical assumption that $M$ be {\sl Levi
nondegenerate} at the point $p$ ({\em see} {\em e.g.} \cite{
isa2009})\,\,---\,\,which is the origin of our present system of
coordinates $(z, w)$\,\,---\,\,may then be expressed here ({\em cf.}
\cite{ me2005a, me2005b}) by requiring that $\Theta_{ \overline{ z}}
\Theta_{ z \overline{ w}} - \Theta_{ \overline{ w}} \Theta_{ z
\overline{ z}}$ does not vanish at the origin. In particular, this
guarantees that the following explicit rational expression whose
numerator is a polynomial in the fourth-order jet $J_{ z, \overline{
z}, \overline{ w}}^4 \Theta$, is well defined and analytic in some
sufficiently small neighborhood of the origin:
\[
\footnotesize
\aligned
{\sf AJ}^4(\Theta)
&
:=
\frac{1}{
[\Theta_{\overline{z}}\Theta_{z\overline{w}}
-\Theta_{\overline{w}}\Theta_{z\overline{z}}]^3}
\bigg\{
\Theta_{zz\overline{z}\overline{z}}
\bigg(
\Theta_{\overline{w}}\Theta_{\overline{w}}
\left\vert\!\!
\begin{array}{cc}
\Theta_{\overline{z}} & \Theta_{\overline{w}}
\\
\Theta_{z\overline{z}} & \Theta_{z\overline{w}}
\end{array}
\!\!\right\vert
\bigg)
-
\\
&
\ \ \ \ \
-\,
2\Theta_{zz\overline{z}\overline{w}}
\bigg(
\Theta_{\overline{z}}\Theta_{\overline{w}}
\left\vert\!\!
\begin{array}{cc}
\Theta_{\overline{z}} & \Theta_{\overline{w}}
\\
\Theta_{z\overline{z}} & \Theta_{z\overline{w}}
\end{array}
\!\!\right\vert
\bigg)
+
\Theta_{zz\overline{w}\overline{w}}
\bigg(
\Theta_{\overline{z}}\Theta_{\overline{z}}
\left\vert\!\!
\begin{array}{cc}
\Theta_{\overline{z}} & \Theta_{\overline{w}}
\\
\Theta_{z\overline{z}} & \Theta_{z\overline{w}}
\end{array}
\!\!\right\vert
\bigg)
+
\\
&
\ \ \ \ \
+
\Theta_{zz\overline{z}}
\bigg(
\Theta_{\overline{z}}\Theta_{\overline{z}}
\left\vert\!\!
\begin{array}{cc}
\Theta_{\overline{w}} & \Theta_{\overline{w}\overline{w}}
\\
\Theta_{z\overline{w}} & \Theta_{z\overline{w}\overline{w}}
\end{array}
\!\!\right\vert
-
2\Theta_{\overline{z}}\Theta_{\overline{w}}
\left\vert\!\!
\begin{array}{cc}
\Theta_{\overline{w}} & \Theta_{\overline{z}\overline{w}}
\\
\Theta_{z\overline{w}} & \Theta_{z\overline{z}\overline{w}}
\end{array}
\!\!\right\vert
+
\Theta_{\overline{w}}\Theta_{\overline{w}}
\left\vert\!\!
\begin{array}{cc}
\Theta_{\overline{w}} & \Theta_{\overline{z}\overline{z}}
\\
\Theta_{z\overline{w}} & \Theta_{z\overline{z}\overline{z}}
\end{array}
\!\!\right\vert
\bigg)
+\\
&
\ \ \ \ \
+
\Theta_{zz\overline{w}}
\bigg(
-\Theta_{\overline{z}}\Theta_{\overline{z}}
\left\vert\!\!
\begin{array}{cc}
\Theta_{\overline{z}} & \Theta_{\overline{w}\overline{w}}
\\
\Theta_{z\overline{z}} & \Theta_{z\overline{w}\overline{w}}
\end{array}
\!\!\right\vert
+
2\Theta_{\overline{z}}\Theta_{\overline{w}}
\left\vert\!\!
\begin{array}{cc}
\Theta_{\overline{z}} & \Theta_{\overline{z}\overline{w}}
\\
\Theta_{z\overline{z}} & \Theta_{z\overline{z}\overline{w}}
\end{array}
\!\!\right\vert
-
\Theta_{\overline{w}}\Theta_{\overline{w}}
\left\vert\!\!
\begin{array}{cc}
\Theta_{\overline{z}} & \Theta_{\overline{z}\overline{z}}
\\
\Theta_{z\overline{z}} & \Theta_{z\overline{z}\overline{z}}
\end{array}
\!\!\right\vert
\bigg)
\bigg\}.
\endaligned
\]
We hope, then, that the following precise statement will fill a gap in our
understanding of the vanishing of CR curvature tensors.

\smallskip\noindent{\bf Main (and unique) theorem.}
{\em 
An arbitrary, not necessarily rigid, real analytic hypersurface $M
\subset \C^2$ which is Levi nondegenerate at one of its points $p$ and
has a complex definining equation of the form:}
\[
w
=
\Theta\big(z,\,\overline{z},\,\overline{w}\big)
\]
{\em in some system of local holomorphic coordinates $(z, w) \in \C^2$
centered at $p$, is spherical at $p$ {\em if and only if} its graphing
complex function $\Theta$ satisfies the following explicit sixth-order
algebraic partial differential equation:}
\[
\small
\aligned
0
\equiv
\bigg(
\frac{-\,\Theta_{\overline{w}}}{
\Theta_{\overline{z}}\Theta_{z\overline{w}}
-\Theta_{\overline{w}}\Theta_{z\overline{z}}}\,
\frac{\partial}{\partial\overline{z}}
+
\frac{\Theta_{\overline{z}}}{
\Theta_{\overline{z}}\Theta_{z\overline{w}}
-\Theta_{\overline{w}}\Theta_{z\overline{z}}}\,
\frac{\partial}{\partial\overline{w}}
\bigg)^2
\big[{\sf AJ}^4(\Theta)\big]
\endaligned
\]
{\em identically in $\C \big\{ z, \overline{ z}, \overline{ w}
\big\}$.}\medskip

Here, it is understood that the first-order derivation in parentheses
is applied twice to the fourth-order rational differential expression
${\sf AJ}^4 ( \Theta)$. The factor $\frac{ 1}{ [ \Theta_{ \overline{
z} } \Theta_{z \overline{ w}} -\Theta_{ \overline{ w } } \Theta_{z
\overline{ z }}]^7}$ then appears, and after clearing out this
denominator, one obtains a universal {\em polynomial} differential
expression ${\sf AJ}^6 ( \Theta)$ depending upon the sixth-order jet
$J_{ z, \overline{ z}, \overline{ w}}^6 \Theta$ and having integer
coefficients. A partial expansion is provided in Section~5, and the
already formidable incompressible length of this expansion perhaps
explains the reason why no reference in the literature provides the
explicit expressions, in terms of some defining function for $M$, of
\'Elie Cartan's two fundamental differential invariants\footnote{\,
{\em See}~\cite{ ca1932} and also~\cite{ ns2003}, where the tight
analogy with second-order ordinary differential equations is well
explained. }
which can (in principle) be used to classify real analytic
hypersurfaces of $\C^2$ up to biholomorphisms, and to at least
characterize sphericality.

Suppose in particular for instance that $M$ is rigid, given by a
complex equation of the form $w = - \overline{ w} + \Xi(z, \overline{
z})$, that is to say with $\Theta ( z, \overline{ z}, \overline{ w})$
of the form $- \overline{ w} + \Xi ( z, \overline{ z})$, so that the
reality condition simply reads here: $\Xi ( z, \overline{ z} ) \equiv
\overline{ \Xi} ( \overline{ z}, z)$. Then as a corollary-exercise,
sphericality is explicitly characterized by a much simpler partial
differential equation that we can write down in expanded form:
\[
\footnotesize
\aligned
0
&
\equiv
\frac{\Xi_{z^2\overline{z}^4}}{
\big(\Xi_{z\overline{z}}\big)^4}
-
6\,
\frac{\Xi_{z^2\overline{z}^3}\,\Xi_{z\overline{z}^2}}{
\big(\Xi_{z\overline{z}}\big)^5}
-
4\,
\frac{\Xi_{z^2\overline{z}^2}\,\Xi_{z\overline{z}^3}}{
\big(\Xi_{z\overline{z}}\big)^5}
-
\frac{\Xi_{z^2\overline{z}}\,\Xi_{z\overline{z}^4}}{
\big(\Xi_{z\overline{z}}\big)^5}
+
\\
&\ \ \ \ \
+
15\,
\frac{\Xi_{z^2\overline{z}^2}\,\big(\Xi_{z\overline{z}^2}\big)^2}{
\big(\Xi_{z\overline{z}}\big)^6}
+
10\,
\frac{\Xi_{z\overline{z}^3}\,\Xi_{z^2\overline{z}}\,\Xi_{z\overline{z}^2}}{
\big(\Xi_{z\overline{z}}\big)^6}
-
15\,
\frac{\Xi_{z^2\overline{z}}\,\big(\Xi_{z\overline{z}^2}\big)^3}{
\big(\Xi_{z\overline{z}}\big)^7},
\endaligned
\]
and this equation should of course hold identically in $\C \big\{ z,
\overline{ z} \big\}$.

\smallskip

Now, here is a summarized description of our arguments of
proof. Beniamino Segre (\cite{ seg1931}) in 1931 and in fact much
earlier Sophus Lie himself in the 1880's ({\em see} {\em
e.g.}~Chapter~10 of Volume~I of the {\em Theorie der
Transformationsgruppen}~\cite{ enlie1888}) showed how to elementarily
associate a unique second-order ordinary differential equation:
\[
w_{zz}(z) 
=
\Phi\big(z,\,w(z),\,w_z(z)\big)
\] 
to the Levi nondegenerate equation $w = \Theta ( z, \overline{ z},
\overline{ w})$ by eliminating the two variables $\overline{ z}$ and
$\overline{ w}$, viewed as parameters, from the two equations $w =
\Theta$ and $w_z = \Theta_z$. We check in great details the
semi-known result that $M$ is spherical at the origin 
if and only if its associated differential equation
is equivalent, under some appropriate local holomorphic point
transformation $(z, w) \longmapsto (z', w') = \big( z' ( z, w), w' (
z, w) \big)$ fixing the origin, to the simplest possible equation $w_{
z' z'} ' ( z' ) = 0$ having null right-hand side, whose obvious
solutions are just the affine complex lines. But since the doctoral
dissertation of Arthur Tresse (defended in 1895 under the direction of
Lie in Leipzig), it is known that, attached to any such differential
equation are two explicit differential invariants:
\[
\aligned
{\sf I}^1
&
:=
\Phi_{w_zw_zw_zw_z}
\ \ \ \ \ \ \ \ \ \
\text{\rm and:}
\\
{\sf I}^2
&
:=
{\sf D}{\sf D}\big(\Phi_{w_zw_z}\big)
-
\Phi_{w_z}\,{\sf D}\big(\Phi_{w_zw_z}\big)
-
4\,{\sf D}\big(\Phi_{ww_z}\big)
+
\\
&
\ \ \ \ \ \ \
+
6\,\Phi_{ww}
-
3\,\Phi_w\,\Phi_{w_zw_z}
+
4\,\Phi_{w_z}\,\Phi_{ww_z},
\\
\text{\rm where}\ \ \
{\sf D}
:=
&
\partial_z
+
w_z\,\partial_w
+
\Phi(z,w,w_z)\,\partial_{w_z},
\endaligned
\]
depending both upon the fourth-order jet of $\Phi$, which, together
with all their covariant differentiations, enable one (in
principle\footnote{\,
To our knowledge, the only existing reference where this strategy is
seriously endeavoured in order to classify second-order ordinary
differential equations $y_{ xx} ( x) = F \big( x, y ( x), y_x ( x)
\big)$ is~\cite{hk1989}, but only for certain point
transformations\,\,---\,\,called there ``{\sl
fiber-preserwing}''\,\,---\,\,of the special form $(x, y) \mapsto (x',
y') = \big( x'(x), \, y'(x, y)\big)$, the first component of which is
independent of $y$.  })
to completely determine when two arbitrarily given differential
equations are equivalent one to another\footnote{\,
Three decades earlier, Christoffel in his famous memoir~\cite{
chri1869} of 1869 devoted to the equivalence problem for Riemannian
metrics discovered that the covariant differentiations of the
curvature provide a full list of differential invariants for positive
definite quadratic infinitesimal metrics.
}. 
A very well-known application is: the vanishing of both ${\sf I}^1$
and ${\sf I}^2$ characterizes equivalence to $w_{ z' z'}' ( z') = 0$.
So in order to characterize sphericality, one only has to reexpress
the vanishing of ${\sf I}^1$ and of ${\sf I}^2$ in terms of the
complex defining function $\Theta ( z, \overline{ z}, \overline{ w})$.
For this, we apply the techniques of computational differential
algebra developed in~\cite{ me2009} which enable us here to explicitly
execute the two-ways transfer between algebraic expressions in the jet
of $\Phi$ and algebraic expressions in the jet of $\Theta$. It then
turns out that the two equations which one obtains by transferring to
$\Theta$ the vanishing of ${\sf I}^1$ and of ${\sf I}^2$ are {\em
conjugate one to another}, so that a single equation suffices, and it
is precisely the one enunciated in the theorem. In fact, this
coincidence is caused by the famous projective duality, explained {\em
e.g.} by Lie and Scheffers in Chapter~10 of~\cite{liesch1893} and
restituted in modern language in~\cite{ brya2000, crsa2005}. It is
indeed well known that to any second-order ordinary differential
equation \thetag{ $\mathcal{ E}$}: $y_{ xx} ( x) = F \big( x, y ( x),
y_x ( x) \big)$ is canonically associated a certain {\sl dual}
second-order ordinary differential equation, call it \thetag{
$\mathcal{ E}^*$}: $b_{ aa} (a) = F^* \big( a, b_a ( a), b_{ aa} ( a)
\big)$, which has the crucial property that:
\[
\aligned
&
{\sf I}_{(\mathcal{E})}^1\ \ \
\text{\small\sf is a nonzero multiple of}\ \ \
{\sf I}_{(\mathcal{E}^*)}^2
\\
\text{\rm and symmetrically also:}\ \ \ \ \ \
&
{\sf I}_{(\mathcal{E})}^2\ \ \
\text{\small\sf is a nonzero multiple of}\ \ \
{\sf I}_{(\mathcal{E}^*)}^1. 
\endaligned
\]
The doctoral dissertation~\cite{ kopp1905} of Koppisch (Leipzig 1905)
cited only {\em passim} by \'Elie Cartan in~\cite{ ca1924} contains
the analytical details of this correspondence, which was well
reconstituted recently in~\cite{ crsa2005} within the context of
projective Cartan connections. But the differential equation which is
dual to the one $w_{ zz} ( z) = \Phi \big( z, \, w ( z), \, w_z ( z)
\big)$ associated to $w = \Theta ( z, \overline{ z}, \overline{ w})$
is easily seen to be just its {\em complex conjugate} \thetag{
$\overline{ \mathcal{ E }}$}: $\overline{ w}_{ \overline{ z} \overline{
z}} (\overline{ z}) = \overline{ \Phi} \big( \overline{ z}, \,
\overline{ w} (\overline{ z}),\, \overline{ w}_{ \overline{ z}}
\big)$, and then as a consequence, ${\sf I}_{ (\overline{ \mathcal{
E}} )}^1 = \overline{ {\sf I}_{ (\mathcal{ E})}^1}$ is the conjugate
of ${\sf I}_{ (\mathcal{ E} )}^1$, and similarly also ${\sf I}_{
(\overline{ \mathcal{ E}} )}^2 = \overline{ {\sf I}_{ (\mathcal{
E})}^2}$ is the conjugate of ${\sf I}_{ (\mathcal{ E} )}^2$. So it is
no mystery that, as said, the sphericality of $M$ at the origin:
\[
0
\equiv
{\sf I}_{(\mathcal{E})}^1
\ \ \ \ \ \
\text{\rm and}
\ \ \ \ \ \
0
\equiv
{\sf I}_{(\mathcal{E})}^2
=
{\small\sf nonzero}\cdot
{\sf I}_{(\overline{\mathcal{E}})}^1
=
{\small\sf nonzero}\cdot
\overline{{\sf I}_{(\mathcal{E})}^1},
\]
can in a simpler way be characterized by the vanishing of the two {\em
mutually conjugate} (complex) equations:
\[
0
\equiv
{\sf I}_{(\mathcal{E})}^1
\ \ \ \ \ \
\text{\rm and}
\ \ \ \ \ \
0
\equiv
\overline{{\sf I}_{(\mathcal{E})}^1},
\] 
which of course amount to just {\em one} (complex) equation. 

\smallskip

To conclude this introduction, we would like to mention firstly that
none of our computations\,\,---\,\,especially in Sections~4
and~5\,\,---\,\,was performed with the help of any computer, and
secondly that the effective characterization of sphericality in higher
complex dimension $n \geqslant 3$ will appear soon~\cite{ me2010b}.

\section*{\S2.~Segre varieties and differential equations}
\label{Section-2}

\subsection*{ Real analytic hypersurfaces in $\C^2$}
Let us consider an arbitrary real analytic hypersurface $M$ in $\C^2$
and let us localize it around one of its points, say $p \in M$. Then
there exist complex affine coordinates:
\[
(z,w)
=
\big(x+iy,\,u+iv\big)
\]
vanishing at $p$ in which $T_p M = \{ u = 0\}$, so that $M$ is
represented in a neighborhood of $p$ by a graphed defining equation of
the form:
\[
u
=
\varphi(x,y,v),
\]
where the real-valued function:
\[
\varphi 
= 
\varphi(x,y,v)
= 
\sum_{k,l,m\in\N\atop
k+l+m\geqslant 2}\,
\varphi_{k,l,m}\,x^k y^lv^m
\in
\R\big\{x,y,u\big\},
\]
which possesses entirely arbitrary real coefficients $\varphi_{
k,l,m}$, vanishes at the origin: $\varphi ( 0) = 0$, together with all
its first order derivatives: $0 = \partial_x \varphi ( 0 ) =
\partial_y \varphi ( 0) = \partial_v \varphi ( 0)$. All studies in the
analytic reflection principle\footnote{\,
The reader might for instance consult the survey~\cite{ mp2006},
pp.~5--44 or the memoirs~\cite{ me2005a, me2005b}, and look also at
some of the concerned references therein. }
show without doubt that the adequate geometric concepts: {\sl Pair of
Segre foliations}, {\sl Segre chains}, {\sl Complexified CR orbits},
{\sl Jets of complexified Segre varietes}, {\sl Rigidity of formal CR
mappings}, {\sl Nondegeneracy conditions}, {\sl CR-reflection
function}, can be viewed correctly only when $M$ is represented by a
so-called {\sl complex defining equation}. Such an equation may be
constructed by simply rewriting the initial real equation of $M$ as:
\[
{\textstyle{\frac{w+\overline{w}}{2}}}
=
\varphi
\big(
{\textstyle{\frac{z+\overline{z}}{2}}},\,
{\textstyle{\frac{z-\overline{z}}{2i}}},\,
{\textstyle{\frac{w-\overline{w}}{2i}}}
\big),
\]
and then by solving\footnote{\,
Thanks to $d \varphi( 0) = 0$, the holomorphic implicit function
theorem readily applies. }
the so written equation with respect to $w$, which yields an equation
of the shape\footnote{\,
Notice that since $d\varphi ( 0) = 0$, one has $\Theta = - \overline{
w} + {\sf order}\,2\,{\sf terms}$. }:
\[
w
=
\Theta\big(z,\,\overline{z},\,\overline{w}\big)
=
\sum_{\alpha,\,\beta,\,\gamma\,\in\,\N\atop
\alpha+\beta+\gamma\geqslant 1}\,
\Theta_{\alpha,\beta,\gamma}\,
z^\alpha\,\overline{z}^\beta\,\overline{w}^\gamma
\in
\C\big\{\overline{z},\,z,\,w\big\},
\]
whose right-hand side converges of course near the origin $(0, 0, 0)
\in \C \times \C \times \C$ and has {\em complex} coefficients
$\Theta_{ \alpha, \beta, \gamma} \in \C$. The paradox that any such
{\em complex} equation provides in fact {\em two} real defining
equations for the {\em real} hypersurface $M$ which is {\em
one}-codimensional, and also in addition the fact that one could as
well have chosen to solve the above equation with respect to
$\overline{ w}$, instead of $w$, these two apparent ``contradictions''
are corrected by means of a fundamental, elementary statement that
transfers to $\Theta$ (in a natural way) the condition of reality:
\[
\overline{\varphi(x,y,u)}
=
\sum_{k+l+m\geqslant 1}\,
\overline{\varphi_{k,l,m}}\,
\overline{x}^k\overline{y}^l\overline{v}^m
=
\sum_{k+l+m\geqslant 1}\,
\varphi_{k,l,m}\,x^ky^lv^m
=
\varphi(x,y,v)
\]
enjoyed by the initial definining function $\varphi$.

\smallskip\noindent{\bf Theorem.}
(\cite{ mp2006}, p.~19\footnote{\,
Compared to~\cite{ mp2006}, we denote here by $\Theta$ the function
denoted there by $\overline{ \Theta}$. })
{\em The complex analytic function $\Theta = \Theta ( z, \overline{
z}, \overline{ w} )$ with $\Theta = - \overline{ w} + {\sf O}(2)$
together with its complex conjugate\footnote{\,
According to a general, common convention, given a power series $\Phi
(t) = \sum_{ \gamma \in \N^n}\, \Phi_\gamma\, t^\gamma$, $t \in \C^n$,
$\Phi_\gamma \in \C$, one defines the series $\overline{ \Phi} (t) :=
\sum_{ \gamma \in \N^n }\, \overline{ \Phi }_\gamma\, t^\gamma$ by
conjugating only its complex coefficients. Then the complex
conjugation operator distributes oneself simultaneously on functions
and on variables: $\overline{ \Phi(t)} \equiv \overline{ \Phi} (\bar
t)$, a trivial property which is nonetheless frequently used in the
formal CR reflection principle (\cite{ me2005a, me2005b}). }:
\[
\overline{\Theta}
=
\overline{\Theta}\big(\overline{z},z,w)
=
\sum_{\alpha,\,\beta,\,\gamma\in\N}\,
\overline{\Theta}_{\alpha,\beta,\gamma}\,
\overline{z}^\alpha\,z^\beta\,w^\gamma
\in
\C\big\{\overline{z},\,z,\,w\big\}
\]
satisfy the two {\em (}equivalent by conjugation{\em )} functional
equations{\em :}}
\begin{equation}
\label{reality-Theta}
\aligned
\overline{w}\equiv 
& \
\overline{\Theta}
\big(\overline{z},z,\Theta(z,\overline{z},\overline{w})\big), 
\\
w\equiv 
& \
\Theta\big(z,\overline{z}, 
\overline{\Theta}(\overline{z},z,w)\big).
\endaligned
\end{equation}
{\em 
Conversely, given a local holomorphic function $\Theta ( z, \overline{
z}, \overline{ w} ) \in \C \{ z, \overline{ z}, \overline{ w} \}$,
$\Theta = - \, \overline{ w} + {\sf O} ( 2)$ which, in conjunction
with its conjugate $\overline{ \Theta} ( \overline{ z}, z, w)$,
satisfies this pair of equivalent identities, then the two zero-sets:}
\[
\big\{
0=-\,w+\Theta\big(z,\,\overline{z},\,\overline{w}\big)
\big\}
\ \ \ \ \ \ \ \ \ \ 
\text{\em and}
\ \ \ \ \ \ \ \ \ \ 
\big\{
0=-\,\overline{w}+
\overline{\Theta}\big(\overline{z},\,z,\,w\big)
\big\}
\] 
{\em coincide and define a local {\em one-codimensional} real analytic
hypersurface $M$ passing through the origin in $\C^2$.}\medskip

As before, let $M$ be an arbitrary real analytic hypersurface passing
through the origin in $\C^2$ equipped with coordinates $(z, w)$, and
assume that $T_0 M = \{ u = 0 \}$. Without loss of generality, we can
and we shall assume that the coordinates are chosen in such a way that
a certain standard convenient normalization condition holds.

\smallskip\noindent{\bf Theorem.}
(\cite{ me2005a}, p.~12) {\em There exists a local complex analytic
change of holomorphic coordinates $h \colon (z, w) \longmapsto (z',
w') = h ( z, w)$ fixing the origin and tangent to the identity at the
origin of the specific form}:
\[
z'
=
z, 
\ \ \ \ \ 
w'
=
g(z,w),
\]
{\em such that the image $M' := h(M)$ has a new complex defining
equation $w' = \Theta' \big( z', \overline{ z}', \overline{ w}' \big)$
satisfying:}
\[
\Theta'\big(0,\overline{z}',\overline{w}'\big)
\equiv 
\Theta'\big(z',0,\overline{w}'\big)
\equiv
-\,\overline{w}',
\]
{\em or equivalently, which has a power series expansion of the form:} 
\[
\small
\aligned
\Theta'\big(z',\overline{z}',\overline{w}'\big)
=
-\,\overline{w}'
+
\sum_{\alpha\geqslant 1,\,\,\beta\geqslant 1}\,
\Theta_{\alpha,\beta,0}'\,
{z'}^\alpha{\overline{z}'}^\beta
+
\sum_{\gamma\geqslant 1}\,{\overline{w}'}^\gamma\,
\sum_{\alpha\geqslant 1,\,\,\beta\geqslant 1}\,
\Theta_{\alpha,\beta,\gamma}'\,
{z'}^\alpha{\overline{z}'}^\beta.
\endaligned
\]

\subsection*{ Levi nondegenerate hypersurfaces}
Leaving aside the real defining equation of $M$, let us now rename the
complex defining equation of $M$ in such normalized coordinates simply
as before: $w = \Theta ( z, \overline{ z}, \overline{ w})$, dropping
all the prime signs. Quite concretely, the real analytic hypersurface
$M$ is said to be {\sl Levi nondegenerate} at the origin if the
coefficient $\Theta_{ 1, 1, 0}$ of $z \overline{ z}$, which may be
checked to always be real because of the reality condition~\thetag{
\ref{reality-Theta}}, is {\em nonzero}. In fact, it is well known that
Levi nondegeneracy is a biholomorphically invariant property, {\em
see} for instance~\cite{ mp2006}, p.~158, but in more conceptual
terms, the following general characterization, which may be taken as a
definition here, holds true. One then readily checks that it is
equivalent to $\Theta_{ 1, 1, 0} \neq 0$ in normalized coordinates.

\smallskip\noindent{\bf Lemma.}
\label{characterization-lndg}
(\cite{ me2005a, me2005b, me2009}) {\em The real analytic hypersurface
$M \subset \C^2$ with $0 \in M$ represented in coordinates $(z, w)$ by
a complex defining equation of the form $w = \Theta ( z, \overline{
z}, \overline{ w})$ is Levi nondegenerate at the origin if and only if
the map:
\[
\big(\overline{z},\overline{w}\big)
\longmapsto
\big(
\Theta(0,\overline{z},\overline{w}),\,\,
\Theta_z(0,\overline{z},\overline{w})
\big)
\]
has nonvanishing $2 \times 2$ Jacobian determinant at $(\overline{ z},
\overline{ w}) = (0, 0)$.
}\medskip

After a possible real dilation of the $z$-coordinate, we can therefore
assume that $\Theta_{ 1, 1, 0} = 1$, and then we are provided with the
following normalization:
\begin{equation}
\label{normal-equation}
w
=
-\,\overline{w}
+
z\overline{z}
+
z\overline{z}\,{\sf O}
\big(\vert z\vert+\vert\overline{w}\vert\big),
\end{equation}
that will be useful shortly. Another, even more convincing argument
for consigning to oblivion the real defining equation $u = \varphi (
x, y, v)$ dates back to Beniamino Segre~\cite{ seg1931}, who observed
that to any real analytic $M$ are associated two deeply linked
objects.

\begin{itemize}

\smallskip\item[{\bf 1)}]
The nowadays so-called {\sl Segre varieties}\footnote{\,
A presentation of the general theory, valuable for generic CR
manifolds of arbitrary codimension $d \geqslant 1$ and of arbitrary CR
dimension $m \geqslant 1$ in $\C^{ m+d}$ enjoying no specific
nondegeneracy condition, may be found in~\cite{ me2005a, me2005b,
mp2006}. }
$S_{ \overline{ q}}$ associated to any point $q \in \C^2$ near the
origin of coordinates $(z_q, w_q)$ that are the complex curves defined
by the equation:
\[
S_{\overline{q}}
:=
\big\{
0
=
-\,w+\Theta\big(z,\,\overline{z}_q,\,\overline{w}_q\big)
\big\},
\]
quite appropriately in terms of the fundamental complex defining
function $\Theta$; this equation is {\em holomorphic} just because its
antiholomorphic terms are set fixed.

\smallskip\item[{\bf 2)}]
When $M$ is Levi nondegenerate at the origin, a second-order {\em
complex} ordinary differential equation\footnote{\,
This idea, usually attributed by contemporary CR geometers to
B.~Segre, dates in fact back (at least) to Chapter~10 of Volume~1 of
the 2\,100 pages long {\em Theorie der Transformationsgruppen} written
by Sophus Lie and Friedrich Engel between 1884 and 1893, where it is
even presented in the uppermost general context. }
of the form:
\[
\label{Segre-Theta-F}
w_{zz}(z)
=
\Phi\big(z,\,w(z),\,w_z(z)\big), 
\]
whose solutions are exactly the Segre varieties of $M$, para\-metrized
by the two initial conditions $w ( 0)$ and $w_z ( 0)$ which correspond
bijectively to the antiholomorphic variables $\overline{ z}_q$ and
$\overline{ w}_q$.

\end{itemize}\smallskip

In fact, the recipe for deriving the second-order differential
equation associated to a local Levi-nondegenerate $M \subset \C^2$
with $0 \in M$ represented by a normalized\footnote{\,
In fact, such a normalization was made in advance just in order to
make things concrete and clear, but thanks to what the Lemma on
p.~\pageref{characterization-lndg} expresses in a biholomorphically
invariant way, everything which follows next holds in an arbitrary
system of coordinates. }
equation of the form~\thetag{ \ref{normal-equation}} is very
simple. Considering that $w = w ( z)$ is given in the equation:
\[
w(z)
=
\Theta\big(z,\,\overline{z},\,\overline{w}\big)
\] 
as a function of $z$ with two supplementary (antiholomorphic)
parameters $\overline{ z}$ and $\overline{ w}$ that one would like to
eliminate, we solve with respect to $\overline{ z}$ and $\overline{
w}$, just by means of the implicit function theorem\footnote{\,
Justification: by our preliminary normalization, the $2 \times 2$
Jacobian determinant $\frac{\partial ( \Theta, \, \Theta_z )}{\partial
( \overline{ z},\, \, \overline{ w})}$ computed at the origin equals
{\scriptsize$\left\vert \!\!\begin{array}{cc} 0 & -1 \\ 1 & 0
\end{array} \!\! \right\vert$}, hence is nonzero.
Without the preliminary normalization, the condition of the Lemma on
p.~\pageref{characterization-lndg} also applies in any case. },
the pair of equations:
\[
\left[
\aligned
w(z)
&
=
\Theta\big(z,\,\overline{z},\,\overline{w}\big)
=
-\,\overline{w}+z\overline{z}
+
z\overline{z}\,{\sf O}\big(\vert z\vert+\vert\overline{w}\vert\big)
\\
w_z(z)
&
=
\Theta_z\big(z,\,\overline{z},\,\overline{w}\big)
=
\overline{z}
+
\overline{z}\,{\sf O}\big(\vert z\vert+\vert\overline{w}\vert\big)
\endaligned\right.
\]
the second one being obtained by differentiating the first one with
respect to $z$, and this yields a representation:
\[
\overline{z}
=
\zeta\big(z,\,w(z),\,w_z(z)\big)
\ \ \ \ \ \ \ \ \ \ \ \ 
\text{\rm and}
\ \ \ \ \ \ \ \ \ \ \ \ 
\overline{w}
=
\xi\big(z,\,w(z),\,w_z(z)\big)
\]
for certain two uniquely defined local complex analytic functions
$\zeta ( z, w, w_z)$ and $\xi ( z, w, w_z)$ of three
complex variables. By means of these functions, we may then replace
$\overline{ z}$ and $\overline{ w}$ in the second derivative:
\[
\aligned
w_{zz}(z)
&
=
\Theta_{zz}
\big(z,\,\overline{z},\,\overline{w}\big)
\\
&
=
\Theta_{zz}
\big(z,\,
\zeta\big(z,\,w(z),\,w_z(z)\big),\,\,
\xi\big(z,\,w(z),\,w_z(z)\big)\big)
\\
&
=:
\Phi\big(z,\,w(z),\,w_z(z)\big),
\endaligned
\]
and this defines without ambiguity the associated differential
equation. More about differential equations will be said in \S3 below.

Of course, any spherical real analytic $M \subset \C^2$ must be Levi
nondegenerate at every point, for the unit $3$-sphere $S^3 \subset
\C^2$ is. It is well known that $S^3$ minus one of its points, for
instance: $S^3 \setminus \{ p_\infty \}$ with $p_\infty := (0,-1)$, is
biholomorphic, through the so-called {\sl Cayley transform}:
\[
(z,w)\longmapsto
\big(
{\textstyle{\frac{i\,z}{1+w}}},\,\,
{\textstyle{\frac{1-w}{2+2\,w}}}
\big)
=:
(z',w')
\ \ \ \ \ 
\text{\rm having inverse:}
\ \ \ \ \ 
(z',w')\longmapsto
\big(
{\textstyle{\frac{-2iz'}{1+2w'}}},\,\,
{\textstyle{\frac{1-2w'}{1+2\,w'}}}
\big)
\]
to the so-called {\sl Heisenberg sphere} of equation:
\[
w'
=
-\,\overline{w}'+z'\overline{z}',
\]
in the target coordinates $(z', w')$, and this model will be more
convenient to deal with for our purposes.

\smallskip\noindent{\bf Proposition.}
\label{characterization-sphericality}
{\em A Levi nondegenerate local real analytic hypersurface $M$ in
$\C^2$ is locally biholomorphic to a piece of the Heisenberg sphere
(hence spherical) if and only if its associated second-order ordinary
complex differential equation is locally equivalent to the Newtonian
free particle equation: $w_{ z'z'} ' = 0$, with identically vanishing
right-hand side.}

\proof 
Indeed, any local equivalence of $M$ to the Heisenberg sphere
transforms its differential equation to the one associated with the
Heisenberg sphere, and then trivially: $w_{ z'}' (z') = \overline{
z}'$, whence $w_{ z'z'} ' (z') = 0$.

Conversely, if the Segre varieties of $M$ are mapped to the solutions
of $w_{ z'z'}' (z') = 0$, namely to the complex affine lines of
$\C^2$, the complex defining equation of the transformed $M'$ must
necessarily be affine:
\begin{equation}
\label{lambda-mu}
w'
= 
\overline{\lambda}'\big(\overline{z}',\,\overline{w}'\big) 
+ 
z'\,
\overline{\mu}'\big(\overline{z}',\,\overline{w}'\big)
=:
\Theta'
\big(z',\overline{z}',\overline{w}'\big),
\end{equation}
with certain coefficients that are holomorphic with respect to
$(\overline{ z}', \overline{ w}')$. Then $\overline{ \lambda}' ( 0) =
0$ since the origin is fixed, and if $\overline{ \mu}' ( 0)$ is
nonzero, one performs the linear transformation $z' \mapsto z'$, $w'
\mapsto w' - \overline{\mu}'(0)\, z'$, which stabilizes both $w_{ z'
z'}' (z') = 0$ and the form of~\thetag{ \ref{lambda-mu}}, to insure
then that $\overline{ \mu}' ( 0 ) = 0$.

Next, the second reality condition~\thetag{ \ref{reality-Theta}} now
reads:
\[
w'
\equiv
\overline{\lambda}'
\big(\overline{z}',\,
\overline{\Theta}'\big(\overline{z}',z',w'\big)\big)
+
z'\,\overline{\mu}'
\big(\overline{z}',\,
\overline{\Theta}'\big(\overline{z}',z',w'\big)\big),
\]
and by differentiating it with respect to $\overline{ z}'$, we get,
without writing the arguments for brevity:
\[
\aligned
0
&
\equiv
\overline{\lambda}_{\overline{z}'}'
+
\overline{\Theta}_{\overline{z}'}'\,
\overline{\lambda}_{\overline{w}'}'
+
z'\,\overline{\mu}_{\overline{z}'}'
+
z'\,\overline{\Theta}_{\overline{z}'}'\,
\overline{\mu}_{\overline{w}'}'
\\
&
\equiv
\overline{\lambda}_{\overline{z}'}'
+
\mu'\,
\overline{\lambda}_{\overline{w}'}'
+
z'\,\overline{\mu}_{\overline{z}'}'
+
z'\,\mu'\,
\overline{\mu}_{\overline{w}'}',
\endaligned
\]
where we replace $\overline{ \Theta}_{ \overline{z}'}'$ in the second
line by its value $\mu' ( z', w')$. But with all the arguments, this
identity reads in full length as the following identity holding in $\C
\big\{ \overline{ z}', z', w' \big\}$:
\[
\aligned
&
-\,\overline{\lambda}_{\overline{z}'}'
\big(\overline{z}',\,
\overline{\Theta}'\big(\overline{z}',z',w'\big)\big) 
-
z'\,\overline{\mu}_{\overline{z}'}'
\big(\overline{z}',\,
\overline{\Theta}'\big(\overline{z}',z',w'\big)\big)
\equiv
\\
&
\ \ \ \ \
\equiv
\mu'(z',w')\,
\overline{\lambda}_{\overline{w}'}'
\big(\overline{z}',\,
\overline{\Theta}'\big(\overline{z}',z',w'\big)\big)
+
z'\,\mu'(z',w')\,
\overline{\mu}_{\overline{w}'}'
\big(\overline{z}',\,
\overline{\Theta}'\big(\overline{z}',z',w'\big)\big).
\endaligned
\]
For convenience, it is better to take $(z', \overline{ z}', \overline{
w}')$ as arguments of this identity instead of $(\overline{ z}', z',
w')$, so we simply replace $w'$ in it by:
\[
\Theta'
\big(
z',\,\overline{z}',\,\overline{w}'
\big), 
\]
we apply the first reality condition~\thetag{ \ref{reality-Theta}} and
we get what we wanted to pursue the reasonings:
\begin{equation}
\label{boxed-identity}
\aligned
-\,\overline{\lambda}_{\overline{z}'}'
\big(\overline{z}',\overline{w}'\big)
-z'\,\overline{\mu}_{\overline{z}'}'
\big(\overline{z}',\overline{w}'\big)
&
\equiv
\mu'\big(
z',\,\overline{\lambda}'(\overline{z}',\overline{w}')
+
z'\,\overline{\mu}'(\overline{z}',\overline{w}')
\big)
\cdot
\\
&
\ \ \ \ \ \ \ \ \ \ \ \ \ \ \ \ \ \ 
\cdot
\Big[
\overline{\lambda}_{\overline{w}'}'
\big(\overline{z}',\overline{w}'\big)
+
z'\,
\overline{\mu}_{\overline{w}'}'
\big(\overline{z}',\overline{w}'\big)
\Big],
\endaligned
\end{equation}
{\em i.e.} an identity holding now in $\C \big\{ z', \overline{z}',
\overline{ w}' \big\}$. The left-hand side being affine with respect
to $z'$, the same must be true of each one of the two factors of the
right-hand side. In particular, the second order derivative of the
first factor with respect to $z'$ must vanish identically:
\[
\aligned
0
&
\equiv
\partial_{z'}\partial_{z'}
\big\{
\mu'\big(z',\overline{\lambda}'+z'\overline{\mu'}\big)
\big\}
\\
&
\equiv
\mu_{z'z'}'
+
2\,\overline{\mu}'\,\mu_{z'w'}'
+
\overline{\mu}'\overline{\mu}'\mu_{w'w'}'. 
\endaligned
\]
Because $M'$ is Levi nondegenerate at the origin, the lemma on
p.~\pageref{characterization-lndg} together with the affine
form~\thetag{ \ref{lambda-mu}} of the defining equation entails that
the map:
\begin{equation}
\label{Jacobian-lambda-mu}
\big(\overline{z}',\overline{w}'\big)
\longmapsto
\big(
\overline{\lambda}'(\overline{z}',\overline{w}'),\,\,
\overline{\mu}'(\overline{z}',\overline{w}')
\big)
\end{equation}
has nonvanishing Jacobian determinant at $(\overline{ z}', \overline{
w}') = (0, 0)$. Consequently, in the above identity
(rewritten with some of the arguments): 
\[
0
\equiv
\mu_{z'z'}'
\big(z',\,\overline{\lambda}'+z'\,\overline{\mu}'\big)
+
2\,\overline{\mu}'\,
\mu_{z'w'}'
\big(z',\,\overline{\lambda}'+z'\,\overline{\mu}'\big)
+
\overline{\mu}'\overline{\mu}'\,
\mu_{w'w'}'
\big(z',\,\overline{\lambda}'+z'\,\overline{\mu}'\big),
\]
we can consider $z'$, $\overline{ \lambda}'$ and $\overline{ \mu }'$
as being just three independent variables. Setting $\overline{ \mu }'
= 0$, we get $0 \equiv \mu_{z' z' }' \big( z',\, \overline{ \lambda }'
\big)$, that is to say: $\mu_{ z' z'} (z', w') \equiv 0$ and then
after division of $\overline{ \mu}'$, we are left with only two terms:
\[
0
\equiv
2\,
\mu_{z'w'}'
\big(z',\,\overline{\lambda}'+z'\,\overline{\mu}'\big)
+
\overline{\mu}'\,
\mu_{w'w'}'
\big(z',\,\overline{\lambda}'+z'\,\overline{\mu}'\big).
\] 
Then again $0 \equiv 2\, \mu_{ z'w'} ( z', w')$ and finally also $0
\equiv \mu_{ w' w'} ( z', w')$. This means that the function:
\[
\mu'(z',w')
=
c_1'z'
+
c_2'w',
\]
with some two constants $c_1', c_2' \in \C$, is {\em linear}.

Now, we claim that $c_2' = 0$ in fact. Indeed, setting $\overline{ z}'
= 0$ in~\thetag{ \ref{boxed-identity}}, we get:
\[
-\,\overline{\lambda}_{\overline{z}'}'
\big(0,\overline{w}'\big)
-
z'\,\overline{c}_1'
\equiv
\big\{
c_1'\,z'+
c_2'\big(
\overline{\lambda}'(0,\overline{w}')+z'\overline{c}_2'\overline{w}'
\big)
\big\}
\cdot
\big[
\overline{\lambda}_{\overline{w}'}'(0,\overline{w}')
+
z'\overline{c}_2'
\big].
\]
The coefficient $c_2' \overline{ c}_2 ' \overline{ c}_2'$ of $(z')^2
\overline{ w}'$ in the right-hand side must vanish, so $c_2' =
0$. Since the rank at the origin of the map~\thetag{
\ref{Jacobian-lambda-mu}} equals $2$, necessarily $\mu' \not\equiv 0$,
so $c_1' \neq 0$, and then $c_1' = 1$ after a suitable dilation of the
$z'$-axis. Next, rewriting the identity~\thetag{
\ref{boxed-identity}}:
\[
-\,\overline{\lambda}_{\overline{z}'}'
\big(\overline{z}',\overline{w}'\big)
-
z'
\equiv
z'\big[\,
\overline{\lambda}_{\overline{w}'}'
\big(\overline{z}',\overline{w}'\big)
\big],
\]
we finally get $\overline{ \lambda}_{ \overline{ z}'}' \equiv
0$ and $\overline{ \lambda}_{\overline{ w}'}' \equiv -1$, 
which means in conclusion that:
\[
\lambda'(z',w')
\equiv
-\,w'
\ \ \ \ \ \ \ \ \ \ \ \ \ \
\text{\rm and}
\ \ \ \ \ \ \ \ \ \ \ \ \ \
\mu'(z',w')
\equiv
z',
\]
so that the equation of $M'$ is the one: $w' = -\, \overline{ w}' +
z'\overline{ z}'$ of the Heisenberg sphere in the target coordinates
$(z', w')$.
\endproof

Thanks to this proposition, in order to characterize the sphericality
of a local real analytic hypersurface $M \subset \C^2$ explicitly in
terms of its complex defining function $\Theta$, our
strategy\footnote{\,
---\,\,indicated already as the accessible Open
Question~2.35 in~\cite{ me2009}\,\,--- 
} 
will be to:

\smallskip$\square$\,\,
characterize the local equivalence to $w_{ z'z'} ' ( z') = 0$ of the
associated differential equation:
\begin{equation}
\label{second-Theta}
w_{zz}(z)
=
\Theta_{zz}
\big(
z,\,
\zeta\big(z,\,w(z),\,w_z(z)\big),\,
\xi\big(z,\,w(z),\,w_z(z)\big)
\big),
\end{equation}
explicitly in terms of the three functions 
$\Theta_{ zz}$, $\zeta$ and $\xi$; 

\smallskip$\square$\,\,
eliminate any occurence of the two auxiliary functions $\zeta$ and
$\xi$ so as to re-express the obtained result only in terms of the
sixth-order jet $J_{ z, \overline{ z}, \overline{ w}}^6
\Theta$. 

\section*{\S3.~Geometry of associated submanifolds of solutions}
\label{Section-3}

The characterization we will obtain holds in fact inside a broader
context than just CR geometry, in terms of what we called in~\cite{
me2009} the {\sl submanifold of solutions} associated to any
second-order ordinary differential equation, no matter whether it
comes or not from a Levi nondegenerate $M \subset \C^2$. In fact, the
elementary foundations towards a general theory embracing all systems
of completely integrable partial differential equations was laid
down~\cite{ me2009}, especially by producing explicit prolongation
formulas for infinitesimal Lie symmetries, with many interesting
problems that are still wide open as soon as the number of
(independent or dependent) variables increases: construction of Cartan
connections; production of differential invariants; full
classification according to the Lie symmetry group.

Fortunately for our present purposes here, the geometry, the
classification, and the Lie transformation group features of second
order ordinary differential equations are essentially completely
understood since the groundbreaking works of Lie~\cite{ lie1883},
followed by a prized thesis by Tresse~\cite{ tr1896} and later by a
celebrated memoir of \'Elie Cartan, {\em see} also~\cite{ gtw1989} and
the references therein.

Accordingly, letting $x \in \K$ and $y \in \K$ be two real or complex
variables (with hence $\K = \R$ or $\C$ throughout), consider any
second-order ordinary differential equation:
\[
y_{xx}(x) 
= 
F\big(x,\,y(x),\,y_x(x)\big) 
\]
having local $\K$-analytic right-hand side $F$, and
denote it by \thetag{
$\mathcal{ E}$} for short. In the space of
first-order jets of arbitrary graphing functions $y = y (x)$ that we
equip with three independent coordinates denoted $(x, y, y_x)$, let us
introduce the vector field:
\[
{\sf D}
:=
\frac{\partial}{\partial x}
+
y_x\,\frac{\partial}{\partial y}
+
F(x,y,y_x)\,\frac{\partial}{\partial y_x},
\]
whose integral curves inside the three-dimensional space $(x, y, y_x)$
correspond, classically, to solving the equation $y_{ xx} (x) = F ( x,
y(x), y_x(x))$ by transforming it into a system of two {\em first-order}
differential equations with the two unknown functions $y (x)$ and $y_x
( x)$.

\smallskip\noindent{\bf Theorem.} 
(\cite{ lie1883, tr1896, ca1924, gtw1989, me2006})
{\em A second-order ordinary differential equation $y_{ xx} = F (x, y,
y_x)$ denoted \thetag{ $\mathcal{ E}$}
with $\K$-analytic right-hand side possesses two fundamental
differential invariants, namely:}
\[
\label{two-invariants}
\aligned
{\sf I}_{(\mathcal{E})}^1
&
:=
F_{y_xy_xy_xy_x}
\ \ \ \ \ \ \ \
\text{\em and:}
\\
{\sf I}_{(\mathcal{E})}^2
&
:=
{\sf D}{\sf D}\big(F_{y_xy_x}\big)
-
F_{y_x}\,{\sf D}\big(F_{y_xy_x}\big)
-
4\,{\sf D}\big(F_{yy_x}\big)
+
\\
&
\ \ \ \ \
+
6\,F_{yy}
-
3\,F_y\,F_{y_xy_x}
+
4\,F_{y_x}\,F_{yy_x}, 
\endaligned
\]
{\em while all other differential invariants are deduced from ${\sf
I}_{ (\mathcal{ E})}^1$ and ${\sf I}_{ (\mathcal{ E})}^2$ by covariant
(in the sense of Tresse) or coframe (in the sense of Cartan)
diffentiations. Moreover, local equivalence to $y_{ x'x'}' (x') = 0$
holds under some invertible local $\K$-analytic point transformation:}
\[
(x,y) 
\longmapsto 
(x',y') 
= 
\big(x'(x,y),\,y'(x,y)\big)
\]
{\em if and only if both invariants vanish:}
\[
0
=
{\sf I}_{(\mathcal{E})}^1
=
{\sf I}_{(\mathcal{E})}^2.
\]

In order to characterize sphericality of an $M \subset \C^2$, it is then
natural and advisable to study what the vanishing of the above two
differential invariants gives when applied to the second order
ordinary differential equation~\thetag{ \ref{second-Theta}} enjoyed by
the defining function $\Theta$. This goal will be pursued
in \S4 below.

\smallskip

For the time being, with the aim of
extending such a kind of characterization to a broader scope, 
following \S2 of~\cite{ me2009}, let us now recall
how one may in a natural way construct a {\sl sumanifold of solutions}
$\mathcal{ M}_{\mathcal{ E}}$ associated to the differential equation
\thetag{ $\mathcal{ E}$} which, when \thetag{ $\mathcal{ E}$} comes
from a Levi nondegenerate local real analytic hypersurface $M \subset
\C^2$, regives without any modification its complex defining equation
$w = \Theta \big( z, \overline{ z}, \overline{ w})$.

To begin with, in the first-order jet space $(x, y, y_x)$ 
that we simply draw as a common three-dimensional space: 

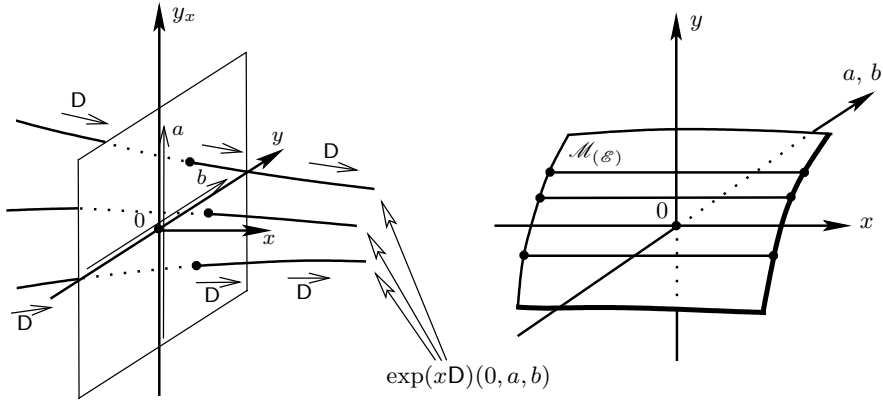
\begin{center}
\begin{picture}(0,0)%
\includegraphics{redressement-A.pstex}%
\end{picture}%
\setlength{\unitlength}{4144sp}%
\begingroup\makeatletter\ifx\SetFigFont\undefined%
\gdef\SetFigFont#1#2#3#4#5{%
  \reset@font\fontsize{#1}{#2pt}%
  \fontfamily{#3}\fontseries{#4}\fontshape{#5}%
  \selectfont}%
\fi\endgroup%
\begin{picture}(5374,2424)(592,-2158)
\put(5720,-1121){\makebox(0,0)[lb]{\smash{{\SetFigFont{9}{10.8}{\familydefault}{\mddefault}{\updefault}{\color[rgb]{0,0,0}$x$}%
}}}}
\put(4501,-1066){\makebox(0,0)[lb]{\smash{{\SetFigFont{9}{10.8}{\familydefault}{\mddefault}{\updefault}{\color[rgb]{0,0,0}$0$}%
}}}}
\put(1753,-862){\makebox(0,0)[lb]{\smash{{\SetFigFont{8}{9.6}{\familydefault}{\mddefault}{\updefault}{\color[rgb]{0,0,0}$b$}%
}}}}
\put(5619,-239){\makebox(0,0)[lb]{\smash{{\SetFigFont{9}{10.8}{\familydefault}{\mddefault}{\updefault}{\color[rgb]{0,0,0}$a,\,b$}%
}}}}
\put(4705, 86){\makebox(0,0)[lb]{\smash{{\SetFigFont{9}{10.8}{\familydefault}{\mddefault}{\updefault}{\color[rgb]{0,0,0}$y$}%
}}}}
\put(2148,-1217){\makebox(0,0)[lb]{\smash{{\SetFigFont{8}{9.6}{\familydefault}{\mddefault}{\updefault}{\color[rgb]{0,0,0}$x$}%
}}}}
\put(2201,-600){\makebox(0,0)[lb]{\smash{{\SetFigFont{8}{9.6}{\familydefault}{\mddefault}{\updefault}{\color[rgb]{0,0,0}$y$}%
}}}}
\put(1605,-568){\makebox(0,0)[lb]{\smash{{\SetFigFont{8}{9.6}{\familydefault}{\mddefault}{\updefault}{\color[rgb]{0,0,0}$a$}%
}}}}
\put(1378,-1117){\makebox(0,0)[lb]{\smash{{\SetFigFont{8}{9.6}{\familydefault}{\mddefault}{\updefault}{\color[rgb]{0,0,0}$0$}%
}}}}
\put(997,-402){\makebox(0,0)[lb]{\smash{{\SetFigFont{8}{9.6}{\familydefault}{\mddefault}{\updefault}{\color[rgb]{0,0,0}${\sf D}$}%
}}}}
\put(678,-1733){\makebox(0,0)[lb]{\smash{{\SetFigFont{8}{9.6}{\familydefault}{\mddefault}{\updefault}{\color[rgb]{0,0,0}${\sf D}$}%
}}}}
\put(1798,-1565){\makebox(0,0)[lb]{\smash{{\SetFigFont{8}{9.6}{\familydefault}{\mddefault}{\updefault}{\color[rgb]{0,0,0}${\sf D}$}%
}}}}
\put(2348,-1546){\makebox(0,0)[lb]{\smash{{\SetFigFont{8}{9.6}{\familydefault}{\mddefault}{\updefault}{\color[rgb]{0,0,0}${\sf D}$}%
}}}}
\put(2509,-702){\makebox(0,0)[lb]{\smash{{\SetFigFont{8}{9.6}{\familydefault}{\mddefault}{\updefault}{\color[rgb]{0,0,0}${\sf D}$}%
}}}}
\put(2889,-2050){\makebox(0,0)[lb]{\smash{{\SetFigFont{9}{10.8}{\familydefault}{\mddefault}{\updefault}{\color[rgb]{0,0,0}$\exp(x{\sf D})(0,a,b)$}%
}}}}
\put(3978,-698){\makebox(0,0)[lb]{\smash{{\SetFigFont{8}{9.6}{\familydefault}{\mddefault}{\updefault}{\color[rgb]{0,0,0}$\mathcal{M}_{(\mathcal{E})}$}%
}}}}
\put(1606,134){\makebox(0,0)[lb]{\smash{{\SetFigFont{9}{10.8}{\familydefault}{\mddefault}{\updefault}{\color[rgb]{0,0,0}$y_x$}%
}}}}
\end{picture}%

\end{center}

\noindent
we {\em duplicate} the two dependent coordinates $(y, y_x)$ by
introducing a new subspace of coordinates $(a, b) \in \K \times \K$,
and we draw a vertical plane containing the two new axes that are
just parallel copies (for the
moment, just look at the left-hand side). 
Then the leaves of the local foliation
associated to the integral curves of the vector field ${\sf D}$ are
uniquely determined by their intersection with this plane, because
thanks to the presence of $\frac{ \partial}{ \partial x}$ in ${\sf
D}$, all these curves are approximately directed by the $x$-axis in a
neighborhood of the origin: no tangent vector can be vertical. But we
claim that all such intersection points of coordinates $(0, b, a) \in
\K \times \K \times \K$ correspond bijectively to the two initial
conditions $y (0) \equiv b$ and $y_x ( 0) = a$ for solving uniquely
the differential equation. In fact, the flow of ${\sf D}$ at time $x$
starting from all such points $(0, b, a)$ of the duplicated vertical
plane:
\[
\exp(x\,{\sf D})(0,b,a)
=:
\big(x,\,Q(x,a,b),\,S(x,a,b)\big)
\]
({\em see} again the diagram) expresses itself
in terms of two certain local
$\K$-analytic functions $Q$ and $S$ that satisfy, by the very
definition of the flow of our vector field $\partial_x + y_x \,
\partial_y + F \, \partial_{ y_x}$, the
following two differential equations:
\[
\frac{d}{dx}
Q(x,a,b)
=
S(x,a,b)
\ \ \ \ \ \ \ \ \ \
\text{\rm and:}
\ \ \ \ \ \ \ \ \ \
\frac{d}{dx}
S(x,a,b)
=
F\big(x,Q(x,a,b),S(x,a,b)\big)
\]
together with the (obious) initial condition for $x = 0$:
\[
(0,b,a)
=
\exp(0\,{\sf D})(0,b,a)
=
\big(0,\,Q(0,a,b),\,S(0,a,b)\big). 
\]
We notice {\em passim} that $S \equiv Q_x$ (no two symbols were
in fact needed), and most importantly, we emphasize that in this way,
we have viewed in a somewhat geometric-minded way of thinking that the
{\em general solution}:
\[
y
=
y(x)
=
Q\big(x,\,y_x(0),\,y(0)\big)
=
Q(x,a,b)
\]
to the original differential equation arises naturally as the
first (amongst two) graphing function for the integral curves of
${\sf D}$ in the first order jet space, these curves being
parametrized by $(a, b)$.

\smallskip\noindent{\bf Definition.}
The {\sl sumanifold of solutions}\footnote{\,
At this point, the reader is referred to~\cite{ me2009} for more about
how one can develope the whole theory of Lie symmetries of partial
differential equations intrinsically within submanifolds of solutions
only; the theory of Cartan connections associated to certain exterior
differential systems could (and should also) be transferred 
to submanifolds of solutions.
} 
$\mathcal{ M}_{(\mathcal{ E})}$ associated with the second-order
ordinary differential equation~\thetag{ $\mathcal{ E}$}:
$y_{ xx} ( x) = F \big( x, \, y (x), \, y_x ( x) \big)$
is the local
$\K$-analytic submanifold of the four-dimensional Euclidean space
$\K_x \times \K_y \times \K_a \times \K_b$ 
represented as the zero-set:
\[
0
=
-y
+
Q(x,a,b), 
\]
where $Q ( x, a,b)$ is the general local $\K$-analytic
solution of~\thetag{ $\mathcal{ E}$}, satisfying therefore: 
\[
Q_{xx}\big(x,a,b)
\equiv
F\big(x,\,Q(x,a,b),\,Q_x(x,a,b)\big), 
\]
and $Q (0, a, b) = b$, $Q_x ( 0, a, b) = a$. 

\medskip
Conversely, let us assume we are given a submanifold $\mathcal{ M}$ of
$\K_x \times \K_y \times \K_a \times \K_b$ of the specific equation $y
= Q ( x, a, b)$, for a certain local $\K$-analytic function $Q$ of
the three variables $(x, a, b)$. Call $(x, y)$ the {\sl variables},
$(a, b)$ the parameters, and call $\mathcal{ M}$ {\sl solvable with
respect to the parameters} (at the origin) if the map:
\[
(a,b)
\longmapsto
\big(Q(0,a,b),\,Q_x(0,a,b)\big)
\] 
has rank two at the central point $(a, b) = (0, 0)$. Of course, the
submanifold of solutions associated to any second-order ordinary
differential equation is solvable with respect to parameters, for in
this case $Q ( 0, a, b) \equiv b$ and $Q_x ( 0, a, b) \equiv a$.

Similarly as what we did for deriving {\bf 2)} on
p.~\pageref{Segre-Theta-F}, if an arbitrarily given submanifold
$\mathcal{ M}$ of $\K_x \times \K_y \times \K_a \times \K_b$ is
assumed to be solvable with respect to parameters, then viewing $y$ in
$y = Q( x, a, b)$ as a parametrized function of $x$, the implicit
function theorem enables one to solve $(a, b)$ in the two equations:
\[
\left[
\aligned
y(x)
&
=
Q(x,a,b)
\\
y_x(x)
&
=
Q_x(x,a,b),
\endaligned\right.
\]
to yield both a representation for $a$ and 
and a representation for $b$ of the form: 
\begin{equation}
\label{A-B}
\left[
\aligned
a
&
=
A\big(x,y(x),y_x(x)\big)
\\
b
&
=
B\big(x,y(x),y_x(x)\big),
\endaligned\right.
\end{equation}
for certain two local $\K$-analytic functions $A$ and $B$ of three
independent variables $(x, y, y_x)$, that one may insert afterwards in
the second order derivative:
\[
\aligned
y_{xx}(x)
&
=
Q_{xx}\big(x,a,b\big)
\\
&
=
Q_{xx}\big(x,\,A(x,y(x),y_x(x)),\,B(x,y(x),y_x(x))\big)
\\
&
=:
F\big(x,\,y(x),\,y_x(x)\big),
\endaligned
\]
which yields the differential equation \thetag{ $\mathcal{
E}_{\mathcal{ M }}$} associated to the submanifold $\mathcal{ M}$
solvable with respect to the parameters. In summary:

\noindent{\bf Proposition.}
(\cite{ me2009})
{\em 
There is a one-to-one correspondence:}
\[
(\mathcal{E}_\mathcal{M})
=
(\mathcal{E})
\longleftrightarrow
\mathcal{M}
=
\mathcal{M}_{(\mathcal{E})}
\]
{\em between second-order ordinary differential equations 
\thetag{ $\mathcal{E}$}
of the
general form:}
\[
y_{xx}(x)
=
F\big(x,\,y(x),\,y_x(x)\big)
\]
{\em and submanifolds {\rm (}of solutions{\rm )} $\mathcal{ M}$ of 
equation:}
\[
y
=
Q(x,a,b)
\]
{\em that are solvable with respect to the parameters, 
and this correspondence satisfies:} 
\[
\big(
\mathcal{E}_{\mathcal{M}_{(\mathcal{E})}}
\big)
=
(\mathcal{E})
\ \ \ \ \ \ \ \ \ \ \ 
\text{\em and}
\ \ \ \ \ \ \ \ \ \ \
\mathcal{M}_{(\mathcal{E}_\mathcal{M})}
=
\mathcal{M}.
\]

We now claim that solvability with respect to the parameters is an
invariant condition, independently of the choice of coordinates.
Indeed, let $y = Q ( x, a, b)$ be any submanifold of solutions, call it
$\mathcal{ M}$, and let: 
\[
\big(x,y,a,b\big)
\longmapsto
\big(x'(x,y),\,y'(x,y),\,a,\,b\big)
\]
be an arbitrary local $\K$-analytic diffeomorphism fixing the origin
which leaves untouched the parameters. The vector of coordinates
$\big(1, \, Q_x ( x, a, b), 0, 0 \big)$ based at the point $\big( x,
Q( x, a, b), \, a, b \big)$ of $\mathcal{ M}$ is sent, through such
a diffeomorphism, to a vector whose $x'$-coordinate equals: $\frac{
d}{ dx} \big[ x' ( x, Q) \big] = x_x' + Q_x \, x_y'$. Therefore
the implicit function theorem insures that, provided the expression:
\[
x_x'(x,y)
+
Q_x(x,a,b)\,x_y'(x,y)
\neq
0
\]
does not vanish, the image $\mathcal{ M}'$ of $\mathcal{ M}$
through such a diffeomorphism can still be represented,
locally in a neighborhood of the origin, 
as a graph of a similar form: 
\[
y'
=
Q'\big(x',\,a,\,b\big),
\]
for a certain local $\K$-analytic new function $Q' = Q' ( x', a,
b)$. Since $\mathcal{ M}: y = Q ( x, a, b)$ is sent to $\mathcal{
M}' : y' = Q' ( x', a, b)$, it follows that $x' (x,y)$, $y' ( x,
y)$, $Q (x, a, b)$ and $Q' (x', a, b)$ are all linked by the
following fundamental identity:
\begin{equation}
y'\big(x,\,Q(x,a,b)\big)
\equiv
Q'\big(x'(x,\,Q(x,a,b)),\,a,b\big),
\end{equation}
which holds in $\C \big\{ x, a, b\big\}$. 

\smallskip\noindent{\bf Claim.}
{\em If $\mathcal{ M}$ is solvable with respect to the parameters (at
the origin), then $\mathcal{ M}'$ is also solvable with respect to the
parameters (at the origin too), and
conversely.}

\proof
The assumption that $\mathcal{ M}$ is solvable with respect to the
parameters is equivalent to the fact that its first order $x$-jet map:
\[
\big(x,\,a,\,b\big)
\longmapsto
\big(
x,\,Q(x,a,b),\,Q_x(x,a,b)\big)
\big)
\]
is (locally) of rank three. One should therefore look at the same
first order jet map attached to $\mathcal{ M}'$, 
represented in the right part of the following diagram:
\[
\footnotesize
\aligned
\diagram 
(x,a,b) \rto 
\dto
&
\big(x'(x,Q(x,a,b)),\,a,\,b\big)
\dto
&
(x',a,b)
\dto
\\
\big(x,\,Q(x,a,b),\,Q_x(x,a,b)\big) \rto^{{\sf X}\,\text{\bf ?}} 
& 
\big(x',\,Q'(x',a,b),\,Q_{x'}'(x',a,b)\big)
&
\!\!\!\!\!\!\!\big(x',\,Q'(x',a,b),\,Q_{x'}'(x',a,b)\big)
\enddiagram,
\endaligned
\]
and ask how these two $x$- and $x'$-jet maps can be related to each
other, namely search for a map: 
\[
{\sf X}\,\text{\bf ?}\colon\,\,
\big(x,\,Q,\,Q_x\big)
\longmapsto
\big(x,\,Q',\,Q_x'\big)
\]
which would close up the diagram and make it commutative.

The answer for the second component of the sought map
is simply: 
\[
{\sf X}_2\colon\,\,
\big(x,\,Q,\,Q_x\big)
\longmapsto
y'\big(x,\,Q\big),
\]
since~\thetag{ 9} indeed shows that composing the right vertical arrow
with the upper horizontal one gives the same result, 
concerning a second
component, as composing the bottom horizontal arrow with the left
vertical one. 

The answer for the third component of the sought map
then proceeds by differentiating with respect to $x$ the fundamental
identity~\thetag{ 9}, which yields, without
writing the arguments: 
\[
y_x'+Q_x\,y_y'
\equiv
\big[x_x'+Q_x\,x_y'\big]\,Q_{x'}', 
\]
and since $x_x' + Q_x \, x_y' \neq 0$ by assumption, 
it suffices to set: 
\[
{\sf X}_3\colon\,\,
\big(x,\,Q,\,Q_x\big)
\longmapsto
\frac{y_x'(x,Q)+Q_x\,y_y'(x,Q)}{x_x'(x,Q)+Q_x\,x_y'(x,Q)},
\]
in order to complete the commutativity of the diagram, 
namely to get:
\[
Q_{x'}'
\big(f(x,Q(x,a,b)),\,a,\,b\big)
\equiv
\frac{y_x'(x,Q(x,a,b))+Q_x(x,a,b)\,y_y'(x,Q(x,a,b))}{
x_x'(x,Q(x,a,b))+Q_x(x,a,b)\,x_y'(x,Q(x,a,b))},
\]
as was required. But now considering instead
the inverse diffeomorphisme changes
nothing to the reasonings, hence we have
at the same time a right-inverse:
\[
\footnotesize
\aligned
\diagram 
(x,a,b) \rto 
\dto^{\text{\rm $x$-jet}}
&
\big(x'(x,Q(x,a,b)),\,a,\,b\big)\rto
\dto^{\text{\rm $x'$-jet}}
&
(x,a,b)
\dto^{\text{\rm $x$-jet}}
\\
\big(x,\,Q(x,a,b),\,Q_x(x,a,b)\big) \rto^{{\sf X}} 
& 
\big(x',\,Q'(x',a,b),\,Q_{x'}'(x',a,b)\big)\rto^{{\sf X}^{-1}}
&
\big(x,\,Q(x,a,b),\,Q(x,a,b)\big)
\enddiagram
\endaligned
\] 
of our commutative diagram, so that the $x$-jet map and the $x'$-jet
map have coinciding ranks at pairs of points which correspond 
one to another. 
\endproof

We are now in a position to generalize the
characterization of sphericality derived
earlier on p.~\pageref{characterization-sphericality}. 

\smallskip\noindent{\bf Proposition.}
{\em A second-order ordinary differential equation $y_{ xx} ( x) = F (
x, y ( x), y_x ( x))$ with $\K$-analytic right-hand side is
equivalent, under some invertible local $\K$-analytic point
transformation $(x, y) \mapsto (x', y')$, to the free particle
Newtonian equation $y_{ x' x'} ' (x') = 0$ if and only if its
associated submanifold of solutions $y = Q( x, a, b)$ is equivalent,
under some local $\K$-analytic map in which variables are separated
from parameters:}
\[
(x,y,a,b)
\longmapsto
\big(x'(x,y),\,y'(x,y),\,\,a'(a,b),\,b'(a,b)\big)
\]
{\em to the affine submanifold of solutions of equation $y' = b' + x'
a'$.}\medskip

Before proceeding to the proof, let us observe that when one looks at
a real analytic hypersurface $M \subset \C^2$, the corresponding
transformation in the parameter space is constrained to be the {\em
conjugate transformation} of the local biholomorphism:
\[
\big(z,w,\,\overline{z},\overline{w}\big)
\longmapsto
\big(
z'(z,w),\,w'(z,w),\,
\overline{z}'(\overline{z},\overline{w}),\,
\overline{w}'(\overline{z},\overline{w})
\big),
\]
while one has more freedom for general differential equations, in the
sense that transformations of variables and transformations of
parameters are entirely {\em decoupled}.

\proof
One direction is clear: if $y = Q ( x, a, b)$ is equivalent to:
\begin{equation}
\label{y-prime}
y'
= 
b'+x'a'
=
b'(a,b)+x'a'(a,b),
\end{equation}
then its associated differential equation $y_{ xx} (x) = F \big( x, y
(x), \, y_x ( x) \big)$ is equivalent, through the same diffeomorphism
$(x, y) \mapsto (x', y')$ of the variables, to the differential
equation associated with~\thetag{ \ref{y-prime}}, which trivially is:
$y_{ x' x'}' ( x') = 0$.

Conversely, if $y_{ xx} (x) = F \big( x, y (x), \, y_x ( x) \big)$ is
equivalent, through a diffeomorphism $(x, y) \mapsto (x', y')$, to
$y_{ x' x'}' ( x') = 0$, then its submanifold of solutions $y = Q (
x, a, b)$ is transformed to $y' = Q' ( x', a, b)$ and since $y_{ x'
x'}' ( x') = 0$, the function $Q'$ is necessarily of the form:
\[
y'
=
b'(a,b)+x'\,a'(a,b). 
\]
Because the condition of solvability with respect to the parameters is
invariant, the rank of $(a, b) \mapsto \big( a' ( a, b), \, b' ( a,
b)\big)$ is again equal to $2$, which concludes the proof.
\endproof

Coming now back to the wanted characterization of sphericality, our
more general goal now amounts to characterize, {\em directly in terms
of its fundamental solution function $Q (x, a, b)$}, the local
equivalence to $y_{ x'x'}' ( x') = 0$ of a second-order ordinary
differential equation $y_{ xx} ( x) = F \big( x, \, y(x), \, y_x ( x)
\big)$. Afterwards at the end, it will suffice to replace $Q ( x, a,
b)$ simply by $\Theta ( z, \overline{ z}, \overline{ w})$ in the
obtained equations.

\smallskip
But before going further, let us explain how a certain generalized
projective duality will simplify our task, as already said in the
Introduction. Thus, let \thetag{ $\mathcal{ E}$}: $y_{ xx} (x) = F
\big( x, y ( x), y_x ( x) \big)$ be a differential equation as above
having general solution $y = Q ( x, a, b) = - b + xa + {\sf O} (
x^2)$, with initial conditions $b = - y ( 0)$ and $a = y_x ( 0)$. The
implicit function theorem enables us to solve $b$ in the equation $y =
Q ( x, a, b)$ of the associated submanifold of solutions $\mathcal{
M}_{ ( \mathcal{ E})}$ in terms of the other quantities, which yields
an equation of the shape:
\[
b
=
Q^*(a,x,y)
=
-y+ax+{\sf O}(x^2), 
\] 
for some new local $\K$-analytic function $Q^* = Q^* ( a, x, y)$. Then
similarly as previously, we may eliminate $x$ and $y$ from the two
equations:
\[
\aligned
b(a)
&
=
Q^*(a,x,y)
=
-y+ax+{\sf O}(x^2)
\\
b_a(a)
&
=
Q_a^*(a,x,y)
=
x+{\sf O}(x^2),
\endaligned
\]
that is to say: $x = X \big( a, b(a), b_a ( a) \big)$ and $y = Y \big(
a, b(a), b_a ( a) \big)$, and we then insert these two
solutions in:
\[
\aligned
b_{aa}(a)
&
=
Q_{aa}^*(a,x,y)
\\
&
=
Q_{aa}^*\big(a,\,X(a,b(a),b_a(a)),\,Y(a,b(a),b_a(a))\big)
\\
&
=:
F^*\big(a,\,b(a),\,b_a(a)\big). 
\endaligned
\]
We shall
call the so obtained second-order ordinary differential equation the
{\sl dual} of $y_{ xx} ( x) = F\big( x, y(x), y_x ( x)\big)$.

In the case of a hypersurface $M \subset \C^2$, solving $\overline{
w}$ in the equation $w = \Theta ( z, \overline{ z}, \overline{ w})$
gives nothing else but the {\em conjugate} equation $\overline{ w} =
\overline{ \Theta} ( \overline{ z}, z, w)$, just by virtue of the reality
identities~\thetag{ \ref{reality-Theta}}. It also follows rather
trivially that the dual differential equation:
\[
\aligned
\overline{w}_{\overline{z}\overline{z}}
(\overline{z})
&
=
\overline{\Theta}_{\overline{z}\overline{z}}
\big(
\overline{z},\,
\overline{\zeta}(\overline{z},\overline{w}(\overline{z}),
\overline{w}_{\overline{z}}(\overline{z})),\,
\overline{\xi}(\overline{z},\overline{w}(\overline{z}),
\overline{w}_{\overline{z}}(\overline{z}))
\big)
\\
&
=
\overline{\Phi}
\big(
\overline{z},\overline{w}(\overline{z}),
\overline{w}_{\overline{z}}(\overline{z})
\big)
\endaligned
\]
is also just the {\em conjugate} differential equation.

To the differential equation $y_{ xx} = F$ and to its dual $b_{ aa} =
F^*$ are associated two submanifolds of solutions:
\[
\mathcal{M}
=
\mathcal{M}_{(\mathcal{E})}
:=
\big\{
\big(x,y,a,b\big)\in\K\times\K\times\K\times\K
\colon\,\,
y=Q(x,a,b)
\big\},
\]
together with: 
\[
\mathcal{M}^*
=
\mathcal{M}_{(\mathcal{E}^*)}
:=
\big\{
\big(a,b,x,y\big)\in\K\times\K\times\K\times\K
\colon\,\,
b=Q^*(a,x,y)
\big\},
\]
and as one obviously guesses, the duality, when viewed within
submanifolds of solutions, just amounts to permute variables and
parameters:
\[
\mathcal{M}
\ni
(x,y,a,b)
\longleftrightarrow
(a,b,x,y)
\in
\mathcal{M}^*.
\]

In the CR case, if we denote by $\widetilde{ z}$ and $\widetilde{ w}$
two independent complex variables which correspond to the
complexifications of $\overline{ z}$ and $\overline{ w}$ (respectively
of course), the duality takes place between the so-called {\sl
extrinsic complexification} (\cite{ me2001b, me2002, me2005a, me2005b,
mp2006, me2009}):
\[
\mathcal{M}
=
M^c
:=
\big\{
\big(z,w,\widetilde{z},\widetilde{w}\big)
\in\C\times\C\times\C\times\C\colon\,\,
w
=
\Theta\big(z,\widetilde{z},\widetilde{w}\big)
\big\}
\]
of $M$ in one hand, and in the other hand, its own
transformation\footnote{\,
Be careful not to write $\big\{ \big( z, w, \widetilde{ z},
\widetilde{ w} \big) \colon \, \, \widetilde{ w} = \overline{ \Theta}
\big( \widetilde{ z}, z, w \big) \big\}$, because this would regive
the same subset $\mathcal{ M}$ of $\C^2 \times \C^2$, due to the
reality identities~\thetag{ \ref{reality-Theta}}.
}: 
\[
\mathcal{M}^*
=
*^c(M^c)
:=
\big\{
\big(\widetilde{z},\widetilde{w},z,w\big)
\in\C\times\C\times\C\times\C\colon\,\,
\widetilde{w}
=
\overline{\Theta}\big(\widetilde{z},z,w\big)
\big\}
\] 
under the involution: 
\[
*^c
\big(z,w,\widetilde{z},\widetilde{w}\big)
:=
\big(\widetilde{z},\widetilde{w},z,w\big)
\]
which clearly is the complexification of the natural antiholomorphic
involution:
\[
*
\big(z,w,\overline{z},\overline{w}\big)
:=
\big(\overline{z},\overline{w},z,w\big)
\]
that fixes $M$ pointwise, as it fixes any other {\em real} analytic
subset of $\C^2$. Here, one has $\mathcal{ M}^* = * ( \mathcal{
M})$\,\,---\,\,which is $\neq \mathcal{ M}$ in general\,\,---\,\,and
of course also $\big( \mathcal{ M}^* \big)^* = \mathcal{ M}$.

So in terms of the coordinates $(x, a, b)$ on $\mathcal{ M}$ and of
the coordinates $(a, x, y)$ on $\mathcal{ M}^*$, the duality is the
map:
\[
(x,a,b)
\longmapsto
\big(a,\,x,\,Q(x,a,b)\big)
\] 
with inverse: 
\[
(a,x,y)
\longmapsto
\big(x,a,Q^*(a,x,y)\big).
\]
But we may also express the duality from the first jet $(x, y,
y_x)$-space to the first jet $(a, b, b_a)$-space by simply composing
the following three maps, the central one being the duality $\mathcal{
M} \to \mathcal{ M}^*$:
\[
\footnotesize
\aligned
\left(
\begin{array}{c}
(a,x,y)
\\
\downarrow
\\
\big(a,Q^*(a,x,y),Q_a^*(a,x,y)\big)
\end{array}
\right)
\circ
\bigg(
(x,a,b)
\rightarrow
\big(a,x,Q(x,a,b)\big)
\bigg)
\circ
\left(
\begin{array}{c}
\big(x,A(x,y,y_x),B(x,y,y_x)\big)
\\
\uparrow
\\
(x,y,y_x)
\end{array}
\right),
\endaligned
\]
which in sum gives us the map: 
\[
(x,y,y_x)
\longmapsto
\left(
\begin{array}{cc}
A(x,y,y_x),\,\,
&
Q^*\big(A(x,y,y_x),x,Q\big(x,A(x,y,y_x),B(x,y,y_x)\big)\big),
\\
&
Q_a^*\big(A(x,y,y_x),x,Q\big(x,A(x,y,y_x),B(x,y,y_x)\big)\big)
\end{array}
\right).
\]
With the approximations, one checks that: 
\[
(x,y,y_x)
\longmapsto
\big(
y_x+\cdots,\,\,
-y+xy_x+\cdots,\,\,
x+\cdots
\big),
\]
where the remainder terms ``$+\cdots$'' are all ${\sf O} (x^2)$. For
the differential equation $y_{ xx} ( x) = 0$ of affine lines, these
remainders disappear completely and we recover the classical
projective duality written in inhomogeneous coordinates (\cite{
crsa2005}, pp.~156--157). Furthermore, one shows ({\em see} {\em
e.g.}~\cite{ crsa2005}) that the above duality map within first order
jet spaces is a {\sl contact transformation}, namely through it, the
pullback of the standard contact form $db - b_a da$ in the target
space is a nonzero multiple of the standard contact form $dy - y_x dx$
in the source space.

But what matters more for us is the following. The two fundamental
differential invariants of $b_{ aa} ( a) = F^* \big( a, b(a), b_a ( a)
\big)$ are functions exactly similar to the ones written on
p.~\pageref{two-invariants}, namely:
\[
\aligned
{\sf I}_{(\mathcal{E}^*)}^1
&
:=
F_{b_ab_ab_ab_a}^*
\\
{\sf I}_{(\mathcal{E}^*)}^2
&
:=
{\sf D}^*{\sf D}^*\big(F_{b_ab_a}^*\big)
-
F_{b_a}^*\,{\sf D}^*\big(F_{b_ab_a}^*\big)
-
4\,{\sf D}^*\big(F_{bb_a}^*\big)
+
\\
&
\ \ \ \ \ \ \ \ \ 
+
6\,F_{bb}^*
-
3\,F_b^*\,F_{b_ab_a}^*
+
4\,F_{b_a}^*\,F_{bb_a}^*,
\endaligned
\]
where ${\sf D}^* := \partial_a + b_a \, \partial_b + F^* ( a, b, b_a)
\, \partial_{ b_a}$. Then according to Koppisch (\cite{ kopp1905}),
through the duality map, ${\sf I}_{ ( \mathcal{ E})}^1$ is transformed to
a nonzero multiple of ${\sf I}_{ ( \mathcal{ E}^*)}^2$, and
simultaneously also, ${\sf I}_{ ( \mathcal{ E})}^2$ is transformed to
a nonzero multiple\footnote{\,
To be precise, both factors of multiplicity (\cite{ crsa2005}, p.~165)
are nonvanishing in a neighborood of the origin, but for our purposes,
it suffices just that they are not identically zero power series.
} 
of ${\sf I}_{ ( \mathcal{ E}^*)}^1$, so that:
\[
\aligned
0={\sf I}_{(\mathcal{E})}^1
\ \ \ \ \
&
\Longleftrightarrow
\ \ \ \ \
{\sf I}_{(\mathcal{E}^*)}^2
=
0
\\
0={\sf I}_{(\mathcal{E})}^2
\ \ \ \ \
&
\Longleftrightarrow
\ \ \ \ \
{\sf I}_{(\mathcal{E}^*)}^1
=
0.
\endaligned
\] 
Consequently, the differential equation \thetag{ $\mathcal{ E }$}:
$y_{ xx} (x) = F \big( x, y ( x), y_x ( x) \big)$ is equivalent to
$y_{ x' x'} ' ( x') = 0$ if and only if:
\[
\boxed{
F_{y_xy_xy_xy_x}
=
0
\ \ \ \ \ \ \ \ \ \
\text{\rm and}
\ \ \ \ \ \ \ \ \ \
F_{b_ab_ab_ab_a}^*
=
0}\,.
\]
This observation has essentially no practical interest, because the
computation of $F^*$ in terms of $F$ relies upon the composition of
three maps \dots\, {\em except notably in the CR case}, since the
duality in this case is complex conjugation: $\Phi^* = \overline{
\Phi}$. In summary, we have established the following.

\smallskip\noindent{\bf Proposition.}
{\em 
An arbitrary, not necessarily rigid, real analytic hypersurface $M
\subset \C^2$ which is Levi nondegenerate at one of its points $p$ and
has a complex definining equation of the form:}
\[
w
=
\Theta\big(z,\,\overline{z},\,\overline{w}\big)
\]
{\em in some system of local holomorphic coordinates $(z, w) \in \C^2$
centered at $p$, is spherical at $p$ {\em if and only if} the
right-hand side $\Phi$ of its uniquely associated second-order
ordinary complex differential equation:}
\[
w_{zz}(z)
=
\Phi\big(z,\,w(z),\,w_z(w)\big)
\]
{\em 
satisfies the {\em single} fourth-order partial differential equation:}
\[
0\equiv
F_{w_zw_zw_zw_z}\big(z,w,w_z\big).
\] 

It now only remains to re-express this fourth-order partial
differential equation in terms of the complex graphing function
$\Theta ( z, \overline{ z}, \overline{ w})$ for $M$.
We will achieve this more generally for $F_{ y_x y_x y_x y_x}$.

\section*{\S4.~Effective differential characterization
\\ 
of sphericality in $\C^2$}
\label{Section-4}

Reminding the reasonings and notations 
introduced in a neighborhood of equation~\thetag{
\ref{A-B}}, the transformation:
\[
\big(x,y,y_x\big)
\longmapsto
\big(x,a,b\big)
\]
and its inverse are given 
by the two triples of functions:
\[
\left[\aligned
x
&
=
x
\\
a
&
=
A(x,y,y_x)
\\
b
&
=
B(x,y,y_x)
\endaligned\right.
\ \ \ \ \ \ \ \ \ \ \ \
\text{\rm and}
\ \ \ \ \ \ \ \ \ \ \ \
\left[\aligned
x
&
=
x
\\
y
&
=
Q(x,a,b)
\\
y_x
&
=
Q_x(x,a,b).
\endaligned\right.
\]
Equivalently, one has the two pairs of identically satisfied
equations:
\[
\small
\aligned
a
&
\equiv
A\big(x,\,Q(x,a,b),\,Q_x(x,a,b)\big)
\\
b
&
\equiv
B\big(x,\,Q(x,a,b),\,Q_x(x,a,b)\big)
\endaligned
\ \ \ \ \
\text{\rm and}
\ \ \ \ \
\aligned
y
&
\equiv
Q\big(x,\,A(x,y,y_x),\,B(x,y,y_x)\big)
\\
y_x
&
\equiv
Q_x\big(x,\,A(x,y,y_x),\,B(x,y,y_x)\big).
\endaligned
\]
Differentiating the second column of equations with respect to $x$,
to $y$ and to $y_x$ yields:
\[
\aligned
0
&
=
Q_x
+
Q_a\,A_x
+
Q_b\,B_x
\ \ \ \ \ \ \ \ \ \ \ \ \ \
0
=
Q_{xx}
+
Q_{xa}\,A_x
+
Q_{xb}\,B_x
\\
1
&
=
\ \ \ \ \ \ \ \ \ \ 
Q_a\,A_y
+
Q_b\,B_y
\ \ \ \ \ \ \ \ \ \ \ \ \ \ 
0
=
\ \ \ \ \ \ \ \ \ \ \
Q_{xa}\,A_y
+
Q_{xb}\,B_y
\\
0
&
=
\ \ \ \ \ \ \ \ \ 
Q_a\,A_{y_x}
+
Q_b\,B_{y_x}
\ \ \ \ \ \ \ \ \ \ \ \ \,
1
=
\ \ \ \ \ \ \ \ \ \ \,
Q_{xa}\,A_{y_x}\!\!
+
Q_{xb}\,B_{y_x}.
\endaligned
\]
Then thanks to a straightforward application
of the rule of Cramer for $2 \times 2$ linear systems,
we derive six useful formulas.

\smallskip\noindent{\bf Lemma.}
(\cite{ me2009}, p.~9)
{\em 
All the six first order derivatives $A_x$, $A_y$, $A_{ y_x}$, $B_x$,
$B_y$, $B_{y_x}$ of the two functions $A$ and $B$ with respect to
their three arguments $(x, y, y_x)$ may be expressed as follows in
terms of the second jet $J^2 (Q)$ of the defining function $Q$:}
\[
\aligned
A_x
&
=
\frac{
Q_b\,Q_{xx}
-
Q_x\,Q_{xb}
}{
Q_a\,Q_{xb}-Q_b\,Q_{xa}},
\ \ \ \ \ \ \ \ \ \ \ \ \ \ \
B_x
=
\frac{
Q_x\,Q_{xa}
-
Q_a\,Q_{xx}
}{
Q_a\,Q_{xb}-Q_b\,Q_{xa}},
\\
A_y
&
=
\frac{
Q_{xb}
}{
Q_a\,Q_{xb}-Q_b\,Q_{xa}},
\ \ \ \ \ \ \ \ \ \ \ \ \ \ \ \,
B_y
=
\frac{
-Q_{xa}
}{
Q_a\,Q_{xb}-Q_b\,Q_{xa}},
\\
A_{y_x}
&
=
\frac{
-Q_b
}{
Q_a\,Q_{xb}-Q_b\,Q_{xa}},
\ \ \ \ \ \ \ \ \ \ \ \ \ \ \,
B_{y_x}
=
\frac{
Q_a
}{
Q_a\,Q_{xb}-Q_b\,Q_{xa}}.
\endaligned
\]

For future abbreviation, we shall denote the single appearing
denominator\,\,---\,\,which evidently is the common 
determinant of all the
three $2 \times 2$ linear systems involved above\,\,---\,\,simply 
by a square symbol:
\[
\Delta
:=
Q_aQ_{xb}
-
Q_bQ_{xa}.
\]
The two-ways transfer between functions $G$ defined in the $(x, y,
y_x)$-space and functions $T$ defined in the $(x, a, b)$-space,
namely the one-to-one correspondence:
\[
G(x,y,y_x)
\longleftrightarrow
T(x,a,b)
\]
may be read very concretely as the following two equivalent identities:
\[
\aligned
G(x,y,y_x)
&
\equiv
T\big(x,\,A(x,y,y_x),\,B(x,y,y_x)\big)
\\
G\big(x,\,Q(x,a,b),\,Q_x(x,a,b)\big)
&
\equiv
T(x,a,b),
\endaligned
\]
holding in $\K \{ x, y, y_x\}$ and in $\K \{ x, a, b\}$ respectively.
By differentiating the first identity, the chain rule shows how
the three first-order derivation operators (basic vector fields)
$\partial_x$, $\partial_y$ and $\partial_{ y_x}$ living in the $(x, y,
y_x)$-space are transformed into the $(x, a, b)$-space:
\[
\footnotesize
\aligned
\frac{\partial}{\partial x}
&
=
\frac{\partial}{\partial x}
+
\bigg(
\frac{Q_b\,Q_{xx}-Q_x\,Q_{xb}}{\Delta}
\bigg)
\frac{\partial}{\partial a}
+
\bigg(
\frac{Q_x\,Q_{xa}-Q_a\,Q_{xx}}{\Delta}
\bigg)
\frac{\partial}{\partial b}
\\
\frac{\partial}{\partial y}
&
=
\ \ \ \ \ \ \ \ \ \ \ \ \ \ \ \ \ \ \ \ \ \ \ \ \ \ \ \ \ \ \ \
\bigg(
\frac{Q_{xb}}{\Delta}
\bigg)
\frac{\partial}{\partial a}
+
\ \ \ \ \ \ \ \ \ \ \ \ \ \ \ \ \ \ \,
\bigg(
\frac{-Q_{xa}}{\Delta}
\bigg)
\frac{\partial}{\partial b}
\\
\frac{\partial}{\partial y_x}
&
=
\ \ \ \ \ \ \ \ \ \ \ \ \ \ \ \ \ \ \ \ \ \ \ \ \ \ \ \ \ \ \,
\bigg(
\frac{-Q_{b}}{\Delta}
\bigg)
\frac{\partial}{\partial a}
+
\ \ \ \ \ \ \ \ \ \ \ \ \ \ \ \ \ \ \ \ \ \ \ 
\bigg(
\frac{Q_{a}}{\Delta}
\bigg)
\frac{\partial}{\partial b}.
\endaligned
\]

\smallskip\noindent{\bf Lemma.}
{\em 
The total differentiation operator ${\sf D} = \partial_x + y_x \,
\partial_y + F \, \partial_{ y_x}$ associated to $y_{ xx} = F ( x,
y,y_x)$ simply transfers to the basic derivation operator along the
$x$-direction:}
\[
{\sf D}
\longleftrightarrow
\partial_x. 
\]

\proof
Reading the three formulas just preceding, by adding the first one to
the second one multiplied by $y_x = Q_x$ together with the third one
multiplied by $F = Q_{ xx}$, one visibly sees that the coefficients
of both $\frac{ \partial}{ \partial a}$ and $\frac{
\partial}{\partial b}$ do vanish in the obtained sum, as announced.
\endproof

Keeping in mind\,\,---\,\,so 
as to avoid any confusion\,\,---\,\,that the same letter
$x$ is used to denote simultaneously the independent variable of the
differential equation $y_{ xx} = F ( x, y, y_x)$ and the non-parameter
variable of the associated submanifold of solutions $y = Q ( x, a,
b)$, we may now write this two-ways transfer ${\sf D}
\longleftrightarrow \partial_x$ exactly as we did in the above three
equations, namely simply as an equality between two derivations living
in the $(x, y, y_x)$-space and in the $(x, a, b)$-space:
\[
{\sf D}
=
\partial_x.
\]

\smallskip\noindent{\bf Lemma.}
{\em 
With $G = G ( x, y, y_x)$ being any local $\K$-analytic function in
the $(x, y, y_x)$-space, the three second-order derivatives $G_{ y_x
y_x}$, $G_{ yy_x}$ and $G_{ y y}$ express as follows in terms of the
second-order jet $J_{ x, a, b}^2 (T)$ of the defining function 
$T$:}
\[
\footnotesize
\aligned
G_{y_xy_x}
&
=
\frac{Q_b\,Q_b}{\Delta^2}\,T_{aa}
-
\frac{2\,Q_a\,Q_b}{\Delta^2}\,T_{ab}
+
\frac{Q_a\,Q_a}{\Delta^2}\,T_{bb}
+
\\
&
\ \ \ \ \
+
\frac{T_a}{\Delta^3}\,
\bigg(
Q_a\,Q_a
\left\vert\!\!
\begin{array}{cc}
Q_b & Q_{bb}
\\
Q_{xb} & Q_{xbb}
\end{array}
\!\!\right\vert
-
2\,Q_a\,Q_b
\left\vert\!\!
\begin{array}{cc}
Q_b & Q_{ab}
\\
Q_{xb} & Q_{xab}
\end{array}
\!\!\right\vert
+
Q_b\,Q_b
\left\vert\!\!
\begin{array}{cc}
Q_b & Q_{aa}
\\
Q_{xb} & Q_{xaa}
\end{array}
\!\!\right\vert
\bigg)
+
\\
&
\ \ \ \ \
+\frac{T_b}{\Delta^3}
\bigg(
-\,Q_a\,Q_a
\left\vert\!\!
\begin{array}{cc}
Q_a & Q_{bb}
\\
Q_{xa} & Q_{xbb}
\end{array}
\!\!\right\vert
+
2\,Q_a\,Q_b
\left\vert\!\!
\begin{array}{cc}
Q_a & Q_{ab}
\\
Q_{xa} & Q_{xab}
\end{array}
\!\!\right\vert
-
Q_b\,Q_b
\left\vert\!\!
\begin{array}{cc}
Q_a & Q_{aa}
\\
Q_{xa} & Q_{xaa}
\end{array}
\!\!\right\vert
\bigg)
\endaligned
\]
\[
\footnotesize
\aligned
G_{yy_x}
&
=
-\,\frac{Q_b\,Q_{xb}}{\Delta^2}\,T_{aa}
+
\frac{Q_a\,Q_{xb}+Q_b\,Q_{xa}}{\Delta^2}\,T_{ab}
-
\frac{Q_a\,Q_{xa}}{\Delta^2}\,T_{bb}
+
\\
&
\ \ \ \ \
+
\frac{T_a}{\Delta^3}
\bigg(
-\,Q_a\,Q_{xa}
\left\vert\!\!
\begin{array}{cc}
Q_b & Q_{bb}
\\
Q_{xb} & Q_{xbb}
\end{array}
\!\!\right\vert
+
\big(Q_a\,Q_{xb}+Q_b\,Q_{xa}\big)
\left\vert\!\!
\begin{array}{cc}
Q_b & Q_{ab}
\\
Q_{xb} & Q_{xab}
\end{array}
\!\!\right\vert
-
Q_b\,Q_{xb}
\left\vert\!\!
\begin{array}{cc}
Q_b & Q_{aa}
\\
Q_{xb} & Q_{xaa}
\end{array}
\!\!\right\vert
\bigg)
+
\\
&
\ \ \ \ \
+\frac{T_b}{\Delta^3}
\bigg(
Q_a\,Q_{xa}
\left\vert\!\!
\begin{array}{cc}
Q_a & Q_{bb}
\\
Q_{xa} & Q_{xbb}
\end{array}
\!\!\right\vert
-
\big(Q_a\,Q_{xb}+Q_b\,Q_{xa}\big)
\left\vert\!\!
\begin{array}{cc}
Q_a & Q_{ab}
\\
Q_{xa} & Q_{xab}
\end{array}
\!\!\right\vert
+
Q_b\,Q_{xb}
\left\vert\!\!
\begin{array}{cc}
Q_a & Q_{aa}
\\
Q_{xa} & Q_{xaa}
\end{array}
\!\!\right\vert
\bigg)
\endaligned
\]
\[
\footnotesize
\aligned
G_{yy}
&
=
\frac{Q_{xb}\,Q_{xb}}{\Delta^2}\,T_{aa}
-
\frac{2\,Q_{xa}\,Q_{xb}}{\Delta^2}\,T_{ab}
+
\frac{Q_{xa}\,Q_{xa}}{\Delta^2}\,T_{bb}
+
\\
&
\ \ \ \ \
+
\frac{T_a}{\Delta^3}
\bigg(
Q_{xa}\,Q_{xa}
\left\vert\!\!
\begin{array}{cc}
Q_b & Q_{bb}
\\
Q_{xb} & Q_{xbb}
\end{array}
\!\!\right\vert
-
2\,Q_{xa}\,Q_{xb}
\left\vert\!\!
\begin{array}{cc}
Q_b & Q_{ab}
\\
Q_{xb} & Q_{xab}
\end{array}
\!\!\right\vert
+
Q_{xb}\,Q_{xb}
\left\vert\!\!
\begin{array}{cc}
Q_b & Q_{aa}
\\
Q_{xb} & Q_{xaa}
\end{array}
\!\!\right\vert
\bigg)
+
\\
&
\ \ \ \ \
+\frac{T_b}{\Delta^3}
\bigg(
-\,Q_{xa}\,Q_{xa}
\left\vert\!\!
\begin{array}{cc}
Q_a & Q_{bb}
\\
Q_{xa} & Q_{xbb}
\end{array}
\!\!\right\vert
+
2\,Q_{xa}\,Q_{xb}
\left\vert\!\!
\begin{array}{cc}
Q_a & Q_{ab}
\\
Q_{xa} & Q_{xab}
\end{array}
\!\!\right\vert
-
Q_{xb}\,Q_{xb}
\left\vert\!\!
\begin{array}{cc}
Q_a & Q_{aa}
\\
Q_{xa} & Q_{xaa}
\end{array}
\!\!\right\vert
\bigg).
\endaligned
\]

\proof
We apply the operator $\frac{ \partial}{ \partial y_x}$, wiewed in the
$(x, a, b)$-space, to the first order derivative $G_{ y_x}$, 
namely we consider:
\[
\partial_{y_x}\big(G_{y_x}\big)
=
\frac{\partial}{\partial y_x}
\bigg[
-\,\frac{Q_b}{\Delta}\,T_a
+
\frac{Q_a}{\Delta}\,T_b
\bigg],
\]
and we then expand carefully the result by collecting 
somewhat in advance the
obtained terms with respect to the derivatives of $T$:
\[
\footnotesize
\aligned
G_{y_xy_x}
&
=
\bigg(
-\,\frac{Q_b}{\Delta}\,
\frac{\partial}{\partial a}
+
\frac{Q_a}{\Delta}\,
\frac{\partial}{\partial b}
\bigg)
\bigg[
-\,\frac{Q_b}{\Delta}\,T_a
+
\frac{Q_a}{\Delta}\,T_b
\bigg]
\\
&
=
\bigg(
\frac{Q_b}{\Delta}\,\frac{Q_{ab}}{\Delta}
-
\frac{Q_b}{\Delta}\,\frac{Q_b\,\Delta_a}{\Delta^2}
\bigg)T_a
+
\frac{Q_b}{\Delta}\,\frac{Q_b}{\Delta}\,T_{aa}
+
\\
&
\ \ \ \ \
+
\bigg(
-\,\frac{Q_b}{\Delta}\,\frac{Q_{aa}}{\Delta}
+
\frac{Q_b}{\Delta}\,\frac{Q_a\,\Delta_a}{\Delta^2}
\bigg)T_b
-
\frac{Q_b}{\Delta}\,\frac{Q_a}{\Delta}\,T_{ab}
+
\\
&
\ \ \ \ \
+
\bigg(
-\,\frac{Q_a}{\Delta}\,\frac{Q_{bb}}{\Delta}
+
\frac{Q_a}{\Delta}\,\frac{Q_b\,\Delta_b}{\Delta^2}
\bigg)T_a
-
\frac{Q_a}{\Delta}\,\frac{Q_b}{\Delta}\,T_{ab}
+
\\
&
\ \ \ \ \
+
\bigg(
\frac{Q_a}{\Delta}\,\frac{Q_{ab}}{\Delta}
-
\frac{Q_a}{\Delta}\,\frac{Q_a\,\Delta_b}{\Delta^2}
\bigg)T_b
+
\frac{Q_a}{\Delta}\,\frac{Q_a}{\Delta}\,
T_{bb}.
\endaligned
\]
The terms involving $T_{ aa}$, $T_{ ab}$, $T_{ bb}$
are exactly the ones exhibited by the lemma for the expression of $G_{
y_x y_x}$. In the four large parentheses which are 
coefficients of $T_a$,
$T_b$, $T_a$, $T_b$, we replace the occurences of
$\Delta_a$, $\Delta_a$, $\Delta_b$, $\Delta_b$ simply by:
\[
\aligned
\Delta_a
&
=
Q_{xb}\,Q_{aa}+Q_a\,Q_{xab}
-
Q_{xa}\,Q_{ab}-Q_b\,Q_{xaa}
\\
\Delta_b
&
=
Q_{xb}\,Q_{ab}+Q_a\,Q_{xbb}
-
Q_{xa}\,Q_{bb}-Q_b\,Q_{xab},
\endaligned
\]
and the total sum of terms coefficiented by $T_a$ in our
expression now becomes:
\[
\footnotesize
\aligned
&
\frac{T_a}{\Delta^3}
\Big(
Q_b\,Q_{ab}
\big[Q_a\,Q_{xb}-Q_b\,Q_{xa}\big]
-
Q_b\,Q_b
\big[
Q_{xb}\,Q_{aa}+Q_a\,Q_{xab}
-
Q_{xa}\,Q_{ab}-Q_b\,Q_{xaa}
\big]
-
\\
&
\ \ \ \ \ \ \
-\,Q_a\,Q_{bb}
\big[Q_a\,Q_{xb}-Q_b\,Q_{xa}\big]
+
Q_a\,Q_b
\big[
Q_{xb}\,Q_{ab}+Q_a\,Q_{xbb}
-
Q_{xa}\,Q_{bb}-Q_b\,Q_{xab}
\big]
\Big)
=
\\
&
=
\frac{T_a}{\Delta^3}
\Big(
Q_a\,Q_b\,Q_{xb}\,Q_{ab}
-
\underline{Q_b\,Q_b\,Q_{xa}\,Q_{ab}}_{
\tiny{\octagon\!\!\!\!\! 1}}
-
Q_b\,Q_b\,Q_{xb}\,Q_{aa}
-
Q_a\,Q_b\,Q_b\,Q_{xab}
+
\\
&
\ \ \ \ \ \ \ \ \ \ \ \ \ \ \ \ \ \ \ \ \ \ \ \ \ \ \ \ \ \ \ \ \ \
\ \ \ \ \ \ \ \ \ \ \ \ \ \ \ \ \ \ \ \ \ \ \ \ \ \ \ \ \ \ \ \ \ \
+
\underline{Q_b\,Q_b\,Q_{xa}\,Q_{ab}}_{
\tiny{\octagon\!\!\!\!\! 1}}
+
Q_b\,Q_b\,Q_b\,Q_{xaa}
-
\\
&
\ \ \ \ \ \ \ \ \ \ \ \ 
-\,Q_a\,Q_a\,Q_{xb}\,Q_{bb}
+
\underline{Q_a\,Q_b\,Q_{xa}\,Q_{bb}}_{
\tiny{\octagon\!\!\!\!\! 2}}
+
Q_a\,Q_b\,Q_{xb}\,Q_{ab}
+
Q_a\,Q_a\,Q_b\,Q_{xbb}
-
\\
&
\ \ \ \ \ \ \ \ \ \ \ \ \ \ \ \ \ \ \ \ \ \ \ \ \ \ \ \ \ \ \ \ \ \
\ \ \ \ \ \ \ \ \ \ \ \ \ \ \ \ \ \ \ \ \ \ \ \ \ \ \ \ \ \ \ \ \ \
-\,\underline{Q_a\,Q_b\,Q_{xa}\,Q_{bb}}_{
\tiny{\octagon\!\!\!\!\! 2}}
-
Q_a\,Q_b\,Q_b\,Q_{xab}
\Big)
=
\\
&
=
\frac{T_a}{\Delta^3}
\Big(
Q_a\,Q_a\big[Q_b\,Q_{xbb}-Q_{xb}\,Q_{bb}\big]
-
2\,Q_a\,Q_b\big[Q_b\,Q_{xab}-Q_{xb}\,Q_{ab}\big]
+
\\
&
\ \ \ \ \ \ \ \ \ \ \ \ \ \ \ \ \ \ \ \ \ \ \ \ \ \ \ \ \ \ \ \ \ \
\ \ \ \ \ \ \ \ \ \ \ \ \ \ \ \ \ \ \ \ \ \ \ \ \ 
+
Q_b\,Q_b\big[Q_b\,Q_{xaa}-Q_{xb}\,Q_{aa}\big]
\Big),
\endaligned
\]
so that we now have effectively reconstituted the three $2 \times 2$
determinants appearing in the second line of the expression claimed by
the lemma for the transfer of $G_{ y_x y_x}$ to the $(x, a,
b)$-space. The treatment of the coefficient of $\frac{ T_b}{
\Delta^3}$ makes only a few differences, hence will be skipped here
(but not in the manuscript). Finally, the two remaining expressions
for $G_{ y y_x}$ and for $G_{ yy}$ are obtained by performing entirely
analogous algebrico-differential computations.
\endproof

\noindent
{\em End of the proof of the Main Theorem.} Applying the above
formula for $G_{ y_x y_x}$ with $x := z$, with $a := \overline{ z}$,
with $b := \overline{ w}$, with $\Delta := \Theta_{ \overline{ z}}
\Theta_{ z \overline{ w}} - \Theta_{ \overline{ w}} \Theta_{ z
\overline{ z}}$, with $G := \Phi$ and with $T := \Theta_{ z z}$, we
exactly get the expression ${\sf AJ}^4 ( \Theta)$ of the Introduction,
and then its further derivative $\partial_{ y_x} \partial_{ y_x} 
\big[ G_{ y_x y_x} \big] = G_{ y_x y_x y_x y_x}$ is exactly:
\[
\footnotesize
\aligned
0
&
\equiv
\bigg(
\frac{-\,\Theta_{\overline{w}}}{
\Theta_{\overline{z}}\Theta_{z\overline{w}}
-\Theta_{\overline{w}}\Theta_{z\overline{z}}}\,
\frac{\partial}{\partial\overline{z}}
+
\frac{\Theta_{\overline{z}}}{
\Theta_{\overline{z}}\Theta_{z\overline{w}}
-\Theta_{\overline{w}}\Theta_{z\overline{z}}}\,
\frac{\partial}{\partial\overline{w}}
\bigg)^2
\big[{\sf AJ}^4(\Theta)\big]
\\
&
=:
\frac{{\sf AJ}^6(\Theta)}{
[\Theta_{\overline{z}}\Theta_{z\overline{w}}
-\Theta_{\overline{w}}\Theta_{z\overline{z}}]^7}. 
\endaligned
\]
As we have said, the vanishing of the second invariant of $w_{ zz} (z)
= \Phi \big( z, w (z), w_z ( z) \big)$ amounts to the complex
conjugation of the above equation, which is then obviously redundant.
Thus, the proof of the Main Theorem is now complete, but we will
nevertheless discuss in a specific final section what ${\sf AJ}^6 (
\Theta)$ would look like in purely expanded form.
\qed

\section*{\S5.~Some complete expansions: 
\\
examples of expression swellings}
\label{Section-5}

Coming back to the non-CR context with the submanifold of solutions
$\mathcal{ M}_{ ( \mathcal{ E})} = \big\{ y = Q ( x, a, b) \big\}$,
let us therefore figure out how to expand the expression
differentiated twice:
\[
\scriptsize
\aligned
G_{y_xy_xy_xy_x}
&
=
\bigg(
-\,\frac{Q_b}{\Delta}\,\frac{\partial}{\partial a}
+
\frac{Q_a}{\Delta}\,\frac{\partial}{\partial b}
\bigg)^2
\bigg\{
\frac{Q_b\,Q_b}{\Delta^2}\,T_{aa}
-
\frac{2\,Q_a\,Q_b}{\Delta^2}\,T_{ab}
+
\frac{Q_a\,Q_a}{\Delta^2}\,T_{bb}
+
\\
&
\ \ \ \ \
+
\frac{T_a}{\Delta^3}\,
\bigg(
Q_a\,Q_a
\left\vert\!\!
\begin{array}{cc}
Q_b & Q_{bb}
\\
Q_{xb} & Q_{xbb}
\end{array}
\!\!\right\vert
-
2\,Q_a\,Q_b
\left\vert\!\!
\begin{array}{cc}
Q_b & Q_{ab}
\\
Q_{xb} & Q_{xab}
\end{array}
\!\!\right\vert
+
Q_b\,Q_b
\left\vert\!\!
\begin{array}{cc}
Q_b & Q_{aa}
\\
Q_{xb} & Q_{xaa}
\end{array}
\!\!\right\vert
\bigg)
+
\\
&
\ \ \ \ \
+\frac{T_b}{\Delta^3}
\bigg(
-\,Q_a\,Q_a
\left\vert\!\!
\begin{array}{cc}
Q_a & Q_{bb}
\\
Q_{xa} & Q_{xbb}
\end{array}
\!\!\right\vert
+
2\,Q_a\,Q_b
\left\vert\!\!
\begin{array}{cc}
Q_a & Q_{ab}
\\
Q_{xa} & Q_{xab}
\end{array}
\!\!\right\vert
-
Q_b\,Q_b
\left\vert\!\!
\begin{array}{cc}
Q_a & Q_{aa}
\\
Q_{xa} & Q_{xaa}
\end{array}
\!\!\right\vert
\bigg)
\bigg\},
\endaligned
\]
which would make the Main Theorem a bit more precise and
explicit.

First of all, we notice that, in the formulas for $G_{ y_x y_x}$, for
$G_{ y y_x}$, for $G_{ yy}$, all the appearing $2 \times 2$
determinants happen to be modifications of the basic Jacobian-like
$\Delta$-determinant:
\[
\Delta\big(a\vert b\big)
:=
\Delta
=
\left\vert\!\!
\begin{array}{cc}
Q_a & Q_b
\\
Q_{xa} & Q_{xb}
\end{array}
\!\!\right\vert,
\]
and we will denote them accordingly by employing the
following (formally and intuitively clear) notations: 
\[
\small
\aligned
&
\Delta\big(b\vert bb\big)
:=
\left\vert\!\!
\begin{array}{cc}
Q_b & Q_{bb}
\\
Q_{xb} & Q_{xbb}
\end{array}
\!\!\right\vert
\ \ \ \ \ \ \ \ \ \ \ \ \
\Delta\big(b\vert ab\big)
:=
\left\vert\!\!
\begin{array}{cc}
Q_b & Q_{ab}
\\
Q_{xb} & Q_{xab}
\end{array}
\!\!\right\vert
\ \ \ \ \ \ \ \ \ \ \ \ \
\Delta\big(b\vert aa\big)
:=
\left\vert\!\!
\begin{array}{cc}
Q_b & Q_{aa}
\\
Q_{xb} & Q_{xaa}
\end{array}
\!\!\right\vert
\\
&
\Delta\big(a\vert bb\big)
:=
\left\vert\!\!
\begin{array}{cc}
Q_a & Q_{bb}
\\
Q_{xa} & Q_{xbb}
\end{array}
\!\!\right\vert
\ \ \ \ \ \ \ \ \ \ \ \ \
\Delta\big(a\vert ab\big)
:=
\left\vert\!\!
\begin{array}{cc}
Q_a & Q_{ab}
\\
Q_{xa} & Q_{xab}
\end{array}
\!\!\right\vert
\ \ \ \ \ \ \ \ \ \ \ \ \
\Delta\big(a\vert aa\big)
:=
\left\vert\!\!
\begin{array}{cc}
Q_a & Q_{aa}
\\
Q_{xa} & Q_{xaa}
\end{array}
\!\!\right\vert,
\endaligned
\]
the bottom line always coinciding with the differentiation with
respect to $x$ of the top line. These abbreviations will be very
appropriate for the next explicit computation, so let us
rewrite the formula for 
$G_{ y_xy_x}$ using this newly introduced
formalism:
\[
\small
\aligned
G_{y_xy_x}
&
=
\frac{1}{\Delta(a\vert b)^3}
\bigg\{
T_{aa}\big[Q_b\,Q_b\,\Delta(a\vert b)\big]
+
T_{ab}\big[-\,2\,Q_a\,Q_b\,\Delta(a\vert b)\big]
+
T_{bb}\big[Q_a\,Q_a\,\Delta(a\vert b)\big]
+
\\
&
\ \ \ \ \ \ \ \ \ \ \ \ \ \ \ \ \ \ \ \ \ \ 
+
T_a
\big[
Q_a\,Q_a\,\Delta(b\vert bb)
-
2\,Q_a\,Q_b\,\Delta(b\vert ab)
+
Q_b\,Q_b\,\Delta(b\vert aa)
\big]
+
\\
&
\ \ \ \ \ \ \ \ \ \ \ \ \ \ \ \ \ \ \ \ \ \ 
+
T_b
\big[
-\,Q_a\,Q_a\,\Delta(a\vert bb)
+
2\,Q_a\,Q_b\,\Delta(a\vert ab)
-
Q_b\,Q_b\,\Delta(a\vert aa)
\big]
\bigg\}.
\endaligned
\]
Then the twelve partial derivatives with respect to $a$ and with
respect to $b$ of all the six determinants $\Delta \big( \! * \!
\vert \! * \! \big)$ appearing in the the second line are easy to
write down:
\[
\small
\aligned
{\textstyle{\frac{\partial}{\partial b}}}
\big[\Delta\big(b\vert bb\big)]
&
=
\underline{\Delta\big(bb\vert bb\big)}_{
\tiny{0\!\!\!\!\!}}
+
\Delta\big(b\vert bbb\big)
\ \ \ \ \ \ \ \ \ \ \ \ \ \ \ \
{\textstyle{\frac{\partial}{\partial a}}}
\big[\Delta\big(b\vert bb\big)]
=
\Delta\big(ab\vert bb\big)
+
\Delta\big(b\vert abb\big)
\\
{\textstyle{\frac{\partial}{\partial b}}}
\big[\Delta\big(b\vert ab\big)]
&
=
\Delta\big(bb\vert ab\big)
+
\Delta\big(b\vert abb\big)
\ \ \ \ \ \ \ \ \ \ \ \ \ \ \ \
{\textstyle{\frac{\partial}{\partial a}}}
\big[\Delta\big(b\vert ab\big)]
=
\underline{\Delta\big(ab\vert ab\big)}_{
\tiny{0\!\!\!\!\!}}
+
\Delta\big(b\vert aab\big)
\\
{\textstyle{\frac{\partial}{\partial b}}}
\big[\Delta\big(b\vert aa\big)]
&
=
\Delta\big(bb\vert aa\big)
+
\Delta\big(b\vert aab\big)
\ \ \ \ \ \ \ \ \ \ \ \ \ \ \ \
{\textstyle{\frac{\partial}{\partial a}}}
\big[\Delta\big(b\vert aa\big)]
=
\Delta\big(ab\vert aa\big)
+
\Delta\big(b\vert aaa\big)
\\
{\textstyle{\frac{\partial}{\partial b}}}
\big[\Delta\big(a\vert bb\big)]
&
=
\Delta\big(ab\vert bb\big)
+
\Delta\big(a\vert bbb\big)
\ \ \ \ \ \ \ \ \ \ \ \ \ \ \ \ \
{\textstyle{\frac{\partial}{\partial a}}}
\big[\Delta\big(a\vert bb\big)]
=
\Delta\big(aa\vert bb\big)
+
\Delta\big(a\vert abb\big)
\\
{\textstyle{\frac{\partial}{\partial b}}}
\big[\Delta\big(a\vert ab\big)]
&
=
\underline{\Delta\big(ab\vert ab\big)}_{
\tiny{0\!\!\!\!\!}}
+
\Delta\big(a\vert abb\big)
\ \ \ \ \ \ \ \ \ \ \ \ \ \ \ \ 
{\textstyle{\frac{\partial}{\partial a}}}
\big[\Delta\big(a\vert ab\big)]
=
\Delta\big(aa\vert ab\big)
+
\Delta\big(a\vert aab\big)
\\
{\textstyle{\frac{\partial}{\partial b}}}
\big[\Delta\big(a\vert aa\big)]
&
=
\Delta\big(ab\vert aa\big)
+
\Delta\big(a\vert aab\big)
\ \ \ \ \ \ \ \ \ \ \ \ \ \ \ \ 
{\textstyle{\frac{\partial}{\partial a}}}
\big[\Delta\big(a\vert aa\big)]
=
\underline{\Delta\big(aa\vert aa\big)}_{
\tiny{0\!\!\!\!\!}}
+
\Delta\big(a\vert aaa\big),
\endaligned
\]
and the underlined terms vanish for the trivial reason that any $2
\times 2$ determinant, two columns of which coincide, vanishes.
Consequently, we may now endeavour the computation of the third order
derivative:
\[
G_{y_xy_xy_x}
=
\bigg(
-\,\frac{Q_b}{\Delta}\,\frac{\partial}{\partial a}
+
\frac{Q_a}{\Delta}\,\frac{\partial}{\partial b}
\bigg)
\big[G_{y_xy_x}\big]. 
\]
When applying the two derivations in parentheses to:
\[
G_{y_xy_x}
=
{\textstyle{\frac{1}{\Delta^3}}}\big\{{\sf expression}\big\}
\]
we start out by differentiating $\frac{ 1}{ \Delta^3}$ multiplied by
${\sf expression}$, and then we differentiate ${\sf expression}$.
Before any contraction, the full expansion of: 
\[
\Delta^5\,G_{y_xy_xy_x}
=
\]
(we indeed clear out the denominator $\Delta^5$) is then:
\[
\tiny
\aligned
&
=
T_{aa}
\big[
3Q_bQ_bQ_b\Delta(a\vert b)\Delta(aa\vert b)
+
3Q_bQ_bQ_b\Delta(a\vert b)\Delta(a\vert ab)
-
3Q_aQ_bQ_b\Delta(a\vert b)\Delta(ab\vert b)
-
3Q_aQ_bQ_b\Delta(a\vert b)\Delta(a\vert bb)
\big]
+
\\
&
+
T_{ab}
\big[
-
6Q_aQ_bQ_b\Delta(a\vert b)\Delta(aa\vert b)
-
6Q_aQ_bQ_b\Delta(a\vert b)\Delta(a\vert ab)
+
6Q_aQ_aQ_b\Delta(a\vert b)\Delta(ab\vert a)
+
6Q_aQ_aQ_b\Delta(a\vert b)\Delta(a\vert bb)
\big]
+
\\
&
+
T_{bb}
\big[
3Q_aQ_aQ_b\Delta(a\vert b)\Delta(aa\vert b)
+
3Q_aQ_aQ_b\Delta(a\vert b)\Delta(a\vert ab)
-
3Q_aQ_aQ_a\Delta(a\vert b)\Delta(ab\vert b)
-
3Q_aQ_aQ_a\Delta(a\vert b)\Delta(a\vert bb)
\big]
+
\endaligned
\]
\[
\tiny
\aligned
+
T_a
&
\Big[
3Q_aQ_aQ_b\Delta(b\vert bb)\Delta(aa\vert b)
+
3Q_aQ_aQ_b\Delta(b\vert bb)\Delta(a\vert ab)
-
3Q_aQ_aQ_a\Delta(b\vert bb)\Delta(ab\vert b)
-
3Q_aQ_aQ_a\Delta(b\vert bb)\Delta(a\vert bb)
-
\\
&
-
6Q_aQ_bQ_b\Delta(b\vert ab)\Delta(aa\vert b)
-
6Q_aQ_bQ_b\Delta(b\vert ab)\Delta(a\vert ab)
+
6Q_aQ_aQ_b\Delta(b\vert ab)\Delta(ab\vert b)
+
6Q_aQ_aQ_b\Delta(b\vert ab)\Delta(a\vert bb)
+
\\
&
+
3Q_bQ_bQ_b\Delta(b\vert aa)\Delta(aa\vert b)
+
3Q_bQ_bQ_b\Delta(b\vert aa)\Delta(a\vert ab)
-
3Q_aQ_bQ_b\Delta(b\vert aa)\Delta(ab\vert b)
-
3Q_aQ_bQ_b\Delta(b\vert aa)\Delta(a\vert bb)
\Big]
+
\endaligned
\]
\[
\tiny
\aligned
+
T_b
&
\Big[
3Q_aQ_aQ_b\Delta(a\vert bb)\Delta(aa\vert b)
+
3Q_aQ_aQ_b\Delta(a\vert bb)\Delta(a\vert ab)
-
3Q_aQ_aQ_a\Delta(a\vert bb)\Delta(ab\vert b)
-
3Q_aQ_aQ_a\Delta(a\vert bb)\Delta(a\vert bb)
-
\\
&
-
6Q_aQ_bQ_b\Delta(a\vert ab)\Delta(aa\vert b)
-
6Q_aQ_bQ_b\Delta(a\vert ab)\Delta(a\vert ab)
+
6Q_aQ_aQ_b\Delta(a\vert ab)\Delta(ab\vert b)
+
6Q_aQ_aQ_b\Delta(a\vert ab)\Delta(a\vert bb)
+
\\
&
+
3Q_bQ_bQ_b\Delta(a\vert aa)\Delta(aa\vert b)
+
3Q_bQ_bQ_b\Delta(a\vert aa)\Delta(a\vert ab)
-
3Q_aQ_bQ_b\Delta(a\vert aa)\Delta(ab\vert b)
-
3Q_aQ_bQ_b\Delta(a\vert aa)\Delta(a\vert bb)
\Big]
+
\endaligned
\]
\[
\tiny
\aligned
+
\Delta(a\vert b)
T_{aaa}
\big[-Q_bQ_bQ_b\Delta(a\vert b)\big]
+
T_{aab}\big[3Q_aQ_bQ_b\Delta(a\vert b)\big]
+
T_{abb}\big[-3Q_aQ_aQ_b\Delta(a\vert b)\big]
+
T_{bbb}\big[Q_aQ_aQ_a\Delta(a\vert b)\big]
+
\endaligned
\]
\[
\tiny
\aligned
+
\Delta(a\vert b)
T_{aa}
\big[
&
-
2Q_bQ_bQ_{ab}\Delta(a\vert b)
-
Q_bQ_bQ_b\Delta(aa\vert b)
-
Q_bQ_bQ_b\Delta(a\vert ab)
+
\\
&
+
2Q_aQ_bQ_b\Delta(a\vert b)
+
Q_aQ_bQ_bQ_b\Delta(ab\vert b)
+
Q_aQ_bQ_b\Delta(a\vert bb)
\big]
+
\endaligned
\]
\[
\tiny
\aligned
+
\Delta(a\vert b)
T_{ab}
\big[
&
2Q_bQ_bQ_{aa}\Delta(a\vert b)
+
2Q_aQ_bQ_{ab}\Delta(a\vert b)
+
2Q_aQ_bQ_b\Delta(aa\vert b)
+
2Q_aQ_bQ_b\Delta(a\vert ab)
-
\\
&
-
2Q_aQ_bQ_{ab}\Delta(a\vert b)
-
2Q_aQ_aQ_{bb}\Delta(a\vert b)
-
2Q_aQ_aQ_b\Delta(ab\vert b)
-
2Q_aQ_aQ_b\Delta(a\vert bb)
\big]
+
\endaligned
\]
\[
\tiny
\aligned
+
\Delta(a\vert b)
T_{bb}
\big[
&
-
2Q_aQ_bQ_{aa}\Delta(a\vert b)
-
Q_aQ_aQ_b\Delta(aa\vert b)
-
Q_aQ_aQ_b\Delta(a\vert ab)
+
\\
&
+
2Q_aQ_aQ_{ab}\Delta(a\vert b)
+
Q_aQ_aQ_a\Delta(ab\vert b)
+
Q_aQ_aQ_a\Delta(a\vert bb)
\big]
+
\endaligned
\]
\[
\tiny
\aligned
+
\Delta(a\vert b)
T_{aa}
\big[
&
-
Q_aQ_aQ_b\Delta(b\vert bb)
+
2Q_aQ_bQ_b\Delta(b\vert ab)
-
Q_bQ_bQ_b\Delta(b\vert aa)
\big]
+
\endaligned
\]
\[
\tiny
\aligned
+
\Delta(a\vert b)
T_{ab}
\big[
&
Q_aQ_aQ_a\Delta(b\vert bb)
-
2Q_aQ_aQ_b\Delta(b\vert ab)
+
Q_aQ_bQ_b\Delta(b\vert aa)
\big]
+
\endaligned
\]
\[
\tiny
\aligned
+
\Delta(a\vert b)
T_{ba}
\big[
&
Q_aQ_aQ_b\Delta(a\vert bb)
-
2Q_aQ_bQ_b\Delta(a\vert ab)
+
Q_bQ_bQ_b\Delta(a\vert aa)
\big]
+
\endaligned
\]
\[
\tiny
\aligned
+
\Delta(a\vert b)
T_{bb}
\big[
&
-Q_aQ_aQ_a\Delta(a\vert bb)
+
2Q_aQ_aQ_b\Delta(a\vert ab)
-
Q_aQ_bQ_b\Delta(a\vert aa)
\big]
+
\endaligned
\]
\[
\tiny
\aligned
+
\Delta(a\vert b)
T_a
\Big[
&
-
2Q_aQ_bQ_{aa}\Delta(b\vert bb)
-
Q_aQ_aQ_b\Delta(ab\vert bb)
-
Q_aQ_aQ_b\Delta(b\vert abb)
+
\\
&
+
2Q_bQ_bQ_{aa}\Delta(b\vert ab)
+
2Q_aQ_bQ_{ab}\Delta(b\vert ab)
+
\underline{2Q_aQ_bQ_b\Delta(ab\vert ab)}_{0}
+
2Q_aQ_bQ_b\Delta(b\vert aab)
-
\\
&
-
2Q_bQ_bQ_{ab}\Delta(b\vert aa)
-
Q_bQ_bQ_b\Delta(ab\vert aa)
-
Q_bQ_bQ_b\Delta(b\vert aaa)
+
\\
&
+
2Q_aQ_aQ_{ab}\Delta(b\vert bb)
+
\underline{Q_aQ_aQ_a\Delta(bb\vert bb)}_{0}
+
Q_aQ_aQ_a\Delta(b\vert bbb)
-
\\
&
-
2Q_aQ_bQ_{ab}\Delta(b\vert ab)
-
2Q_aQ_aQ_{bb}\Delta(b\vert ab)
-
2Q_aQ_aQ_b\Delta(bb\vert ab)
-
2Q_aQ_aQ_b\Delta(b\vert abb)
+
\\
&
+
2Q_aQ_bQ_{bb}\Delta(b\vert aa)
+
Q_aQ_bQ_b\Delta(bb\vert aa)
+
Q_aQ_bQ_b\Delta(b\vert aab)
\Big]
+
\endaligned
\]
\[
\tiny
\aligned
+
\Delta(a\vert b)
T_b
\Big[
&
2Q_aQ_bQ_{aa}\Delta(a\vert bb)
+
Q_aQ_aQ_b\Delta(aa\vert bb)
+
Q_aQ_aQ_b\Delta(a\vert abb)
-
\\
&
-
2Q_bQ_bQ_{aa}\Delta(a\vert ab)
-
2Q_aQ_bQ_{ab}\Delta(a\vert ab)
-
\underline{2Q_aQ_bQ_b\Delta(aa\vert ab)}_{0}
-
2Q_aQ_bQ_b\Delta(a\vert aab)
+
\\
&
+
2Q_bQ_bQ_{ab}\Delta(a\vert aa)
+
Q_bQ_bQ_b\Delta(aa\vert aa)
+
Q_bQ_bQ_b\Delta(a\vert aaa)
-
\\
&
-
2Q_aQ_aQ_{ab}\Delta(a\vert bb)
-
\underline{Q_aQ_aQ_a\Delta(ab\vert bb)}_{0}
-
Q_aQ_aQ_a\Delta(a\vert bbb)
+
\\
&
+
2Q_aQ_bQ_{ab}\Delta(a\vert ab)
+
2Q_aQ_aQ_{bb}\Delta(a\vert ab)
+
2Q_aQ_aQ_b\Delta(ab\vert ab)
+
2Q_aQ_aQ_b\Delta(a\vert abb)
-
\\
&
-
2Q_aQ_bQ_{bb}\Delta(a\vert aa)
-
Q_aQ_bQ_b\Delta(ab\vert aa)
-
Q_aQ_bQ_b\Delta(a\vert aab)
\Big].
\endaligned
\]
The simplification (collecting all terms) gives: 
\[
\scriptsize
\aligned
G_{y_xy_xy_x}
=
\frac{1}{[\Delta(a\vert b)]^5}
\bigg\{
\aligned
&\ \ \
T_{aaa}
\big[
-
Q_b^3\Delta(a\vert b)^2
\big]
+
T_{aab}
\big[
3Q_aQ_b^2\Delta(a\vert b)^2
\big]
+
\\
&
+
T_{abb}
\big[
-
3Q_a^2Q_b\Delta(a\vert b)^2
\big]
+
T_{bbb}
\big[
Q_a^3\Delta(a\vert b)^2
\big]
+
\endaligned
\endaligned
\]
\[
\scriptsize
\aligned
+
T_{aa}
\Big[
&
-
2Q_b^2Q_{ab}\Delta(a\vert b)^2
+
2Q_aQ_bQ_{bb}\Delta(a\vert b)^2
+
3Q_b^3\Delta(a\vert b)\Delta(aa\vert b)
+
2Q_b^3\Delta(a\vert b)\Delta(a\vert ab)
-
\\
&
-
4Q_aQ_b^2\Delta(a\vert b)\Delta(ab\vert b)
-
2Q_aQ_b^2\Delta(a\vert b)\Delta(a\vert bb)
-
Q_a^2Q_b\Delta(a\vert b)\Delta(b\vert bb)
\Big]
+
\endaligned
\]
\[
\scriptsize
\aligned
+
T_{ab}
\Big[
&
-
2Q_a^2Q_{bb}\Delta(a\vert b)^2
+
2Q_bQ_bQ_{aa}\Delta(a\vert b)^2
+
Q_a^3\Delta(a\vert b)\Delta(b\vert bb)
+
6Q_a^2Q_b\Delta(a\vert b)\Delta(ab\vert b)
+
\\
&
+
Q_b^3\Delta(a\vert b)\Delta(a\vert aa)
-
6Q_aQ_b^2\Delta(a\vert b)\Delta(a\vert ab)
+
5Q_a^2Q_b\Delta(a\vert b)\Delta(a\vert bb)
-
5Q_aQ_b^2\Delta(a\vert b)\Delta(aa\vert b)
\Big]
+
\endaligned
\]
\[
\scriptsize
\aligned
+
T_{bb}
\Big[
&
-
2Q_aQ_bQ_{aa}\Delta(a\vert b)^2
+
2Q_a^2Q_{ab}\Delta(a\vert b)^2
-
3Q_a^3\Delta(a\vert b)\Delta(a\vert bb)
-
2Q_a^3\Delta(a\vert b)\Delta(ab\vert b)
+
\\
&
+
4Q_a^2Q_b\Delta(a\vert b)\Delta(a\vert ab)
+
2Q_a^2Q_b\Delta(a\vert b)\Delta(aa\vert b)
-
Q_aQ_b^2\Delta(a\vert b)\Delta(a\vert aa)
\Big]
+
\endaligned
\]
\[
\tiny
\aligned
+
T_a
\Big[
&
3Q_a^2Q_b\Delta(aa\vert b)\Delta(b\vert bb)
+
3Q_a^2Q_b\Delta(a\vert ab)\Delta(b\vert bb)
-
3Q_a^3\Delta(ab\vert b)\Delta(b\vert bb)
-
3Q_a^3\Delta(a\vert bb)\Delta(b\vert bb)
-
\\
&
-
6Q_aQ_b^2\Delta(aa\vert b)\Delta(b\vert ab)
-
6Q_aQ_b^2\Delta(a\vert ab)\Delta(b\vert ab)
+
6Q_a^2Q_b\Delta(ab\vert b)\Delta(b\vert ab)
+
6Q_a^2Q_b\Delta(a\vert bb)\Delta(b\vert ab)
-
\\
&
-
3Q_b^3\Delta(aa\vert b)\Delta(b\vert aa)
-
3Q_b^3\Delta(a\vert ab)\Delta(b\vert aa)
+
3Q_aQ_b^2\Delta(ab\vert b)\Delta(b\vert aa)
+
3Q_aQ_b^2\Delta(a\vert bb)\Delta(b\vert aa)
-
\\
&
-
2Q_aQ_bQ_{aa}\Delta(a\vert b)\Delta(b\vert bb)
+
2Q_b^2Q_{aa}\Delta(a\vert b)\Delta(b\vert ab)
+
2Q_aQ_bQ_{ab}\Delta(a\vert b)\Delta(b\vert ab)
-
2Q_b^2Q_{ab}\Delta(a\vert b)\Delta(b\vert aa)
-
\\
&
-
Q_a^2Q_b\Delta(a\vert b)\Delta(ab\vert bb)
-
Q_a^2Q_b\Delta(a\vert b)\Delta(b\vert abb)
+
2Q_aQ_b^2\Delta(a\vert b)\Delta(b\vert aab)
-
Q_b^3\Delta(a\vert b)\Delta(ab\vert aa)
-
Q_b^3\Delta(a\vert b)\Delta(b\vert aaa)
+
\\
&
+
2Q_a^2Q_{ab}\Delta(a\vert b)\Delta(b\vert bb)
-
2Q_aQ_bQ_{ab}\Delta(a\vert b)\Delta(b\vert ab)
-
2Q_a^2Q_{bb}\Delta(a\vert b)\Delta(b\vert ab)
+
2Q_aQ_bQ_{bb}\Delta(a\vert b)\Delta(b\vert aa)
+
\\
&
+
Q_a^3\Delta(a\vert b)\Delta(b\vert bbb)
-
2Q_a^2Q_b\Delta(a\vert b)\Delta(bb\vert ab)
-
2Q_a^2Q_b\Delta(a\vert b)\Delta(b\vert abb)
+
Q_aQ_b^2\Delta(a\vert b)\Delta(bb\vert aa)
+
Q_aQ_b^2\Delta(a\vert b)\Delta(b\vert aab)
\Big]
+
\endaligned
\]
\[
\tiny
\aligned
+
T_b
\Big[
&
3Q_a^2Q_b\Delta(a\vert bb)\Delta(aa\vert b)
+
3Q_a^2Q_b\Delta(a\vert bb)\Delta(a\vert ab)
-
3Q_a^3\Delta(a\vert bb)\Delta(ab\vert b)
-
3Q_a^3\Delta(a\vert bb)\Delta(a\vert bb)
-
\\
&
-
6Q_aQ_b^2\Delta(a\vert ab)\Delta(aa\vert b)
-
6Q_aQ_b^2\Delta(a\vert ab)\Delta(a\vert ab)
+
6Q_a^2Q_b\Delta(a\vert ab)\Delta(ab\vert b)
+
6Q_aQ_aQ_b\Delta(a\vert ab)\Delta(a\vert bb)
+
\\
&
+
3Q_b^2\Delta(a\vert aa)\Delta(aa\vert b)
+
3Q_b^2\Delta(a\vert aa)\Delta(a\vert ab)
-
3Q_aQ_b^2\Delta(a\vert aa)\Delta(ab\vert b)
-
3Q_aQ_b^2\Delta(a\vert aa)\Delta(a\vert bb)
+
\\
&
+
2Q_aQ_bQ_{aa}\Delta(a\vert b)\Delta(a\vert bb)
-
2Q_b^2Q_{aa}\Delta(a\vert b)\Delta(a\vert ab)
-
2Q_aQ_bQ_{ab}\Delta(a\vert b)\Delta(a\vert ab)
+
2Q_b^2Q_{ab}\Delta(a\vert b)\Delta(a\vert aa)
+
\\
&
+
Q_a^2Q_b\Delta(a\vert b)\Delta(aa\vert bb)
+
Q_a^2Q_b\Delta(a\vert b)\Delta(a\vert abb)
-
2Q_aQ_b^2\Delta(a\vert b)\Delta(a\vert aab)
+
Q_b^3\Delta(a\vert b)\Delta(aa\vert aa)
+
Q_b^3\Delta(a\vert b)\Delta(a\vert aaa)
-
\\
&
-
2Q_a^2Q_{ab}\Delta(a\vert b)\Delta(a\vert bb)
+
2Q_aQ_bQ_{ab}\Delta(a\vert b)\Delta(a\vert ab)
+
2Q_a^2Q_{bb}\Delta(a\vert b)\Delta(a\vert ab)
-
2Q_aQ_bQ_{bb}\Delta(a\vert b)\Delta(a\vert aa)
-
\\
&
-
Q_a^3\Delta(a\vert b)\Delta(a\vert bbb)
+
2Q_a^2Q_b\Delta(a\vert b)\Delta(ab\vert ab)
+
2Q_a^2Q_b\Delta(a\vert b)\Delta(a\vert abb)
-
Q_aQ_b^2\Delta(a\vert b)\Delta(ab\vert aa)
-
Q_aQ_b^2\Delta(a\vert b)\Delta(a\vert aab)
\Big].
\endaligned
\]
The full expansion of $G_{ y_xy_xy_xy_x}$ will not
be presented here.

\vfill\end{document}